\documentclass[12pt]{article}
\usepackage{hyperref}

\usepackage[utf8]{inputenc}
\usepackage[T1]{fontenc}

\usepackage[]{amsfonts}
\usepackage{amssymb}
\usepackage{amsmath}
\def\couleur(#1 #2 #3)
	{
\special{ps::[end]}}

\def\bx#1{\setbox1=\hbox{\kern3pt{#1}\kern3pt}			
 \dimen1=\ht1 \advance\dimen1 by 3pt \dimen2=\dp1 \advance\dimen2 by 3pt
 \setbox1=\hbox{\vrule height\dimen1 depth\dimen2\box1\vrule}%
 \setbox1=\vbox{\hrule\box1\hrule}%
 \advance\dimen1 by .4pt \ht1=\dimen1
 \advance\dimen2 by .4pt \dp1=\dimen2 \box1\relax}

\def\wbb#1{\kern#1em}
\def\vci{\vrule  width.02em height1.47ex depth-.0ex}		
\def\11{{\rm\wbb{.2}\vci\wbb{-.37}1}}

\def\underset#1#2{\mathrel{\mathop{\kern0pt #2}\limits_{#1}}}

\def\overset#1#2{\mathrel{\mathop{\kern0pt #2}\limits^{#1}}}

\newtheorem{Thrm}{Theorem}[section]
\newtheorem{Lmm}[Thrm]{Lemma}
\newtheorem{Dfnt}[Thrm]{Definition}
\newtheorem{Prps}[Thrm]{Proposition}
\newtheorem{Crll}[Thrm]{Corollary}
\newtheorem{Rmrq}[Thrm]{Remark}

\begin{document}

\title{Nevanlinna classes associated to a closed set on $\displaystyle \partial {\mathbb{D}}.$}

\author{Eric Amar}

\date{ }
\maketitle
 \ \par 
\ \par 
\renewcommand{\abstractname}{Abstract}

\begin{abstract}
We introduce Nevanlinna classes of holomorphic functions associated
 to a closed set on the boundary of the unit disc in the complex
 plane and we get Blaschke type theorems relative to these classes
 by use of  several complex variables methods. This gives alternative
 proofs of some results of Favorov \&  Golinskii, useful, in
 particular, for the study of eigenvalues of non self adjoint
 Schr\"odinger operators.\ \par 
\end{abstract}

\tableofcontents
\ \par 

\section{Introduction.}
\quad We shall study classes of holomorphic functions whose zeros may
 appear as eigenvalues of a Schr\"odinger operator with a non
 self adjoint potential. For instance Frank and Sabin~\cite{FranSabi14}
 use the work of Boritchev, Golinskii and Kupin~\cite{BGK09}
 to get interesting estimates this way.\ \par 
\quad The aim of this work is to study Blaschke type conditions relative
 to Nevanlinna classes associated to a closed set on the torus.
 In order to do this we shall use the "way of thinking of several
 complex variables".\ \par 
\quad The methods used in several complex variables already proved
 their usefulness in the one variable case. For instance:\ \par 
$\bullet $ the corona theorem of Carleson~\cite{CarlesonCor}
 is easier to prove and to understand thanks to the proof of
 T. Wolff based on L. H\"ormander~\cite{HormCor67} ;\ \par 
$\displaystyle \bullet $ the characterization of interpolating
 sequences by Carleson for $\displaystyle H^{\infty }$  and by
 Shapiro \&  Shields for $\displaystyle H^{p}$ are also easier
 to prove by these methods (see~\cite{amExt83}, last section,
 where they allow me to get the bounded linear extension property
 for the case $\displaystyle H^{p}$ ; the $\displaystyle H^{\infty
 }$ case being done by Pehr Beurling~\cite{PBeurling62}).\ \par 
\quad So it is not too surprising that in the case of zero sets, they
 can also be useful.\ \par 
\ \par 
\quad In this work we shall define Nevanlinna classes of holomorphic
 functions in the unit disc $\displaystyle {\mathbb{D}}$ of $\displaystyle
 {\mathbb{C}}$ associated to a closed set $\displaystyle E$ in
 the torus $\displaystyle {\mathbb{T}}$ and we show that the
 zero set of  functions in these Nevanlinna classes must satisfy
 a Blaschke type condition.\ \par 
\quad In fact, the only thing we use with respect to $\displaystyle
 u=\log \left\vert{f(z)}\right\vert $ is the fact that $u$ is
 a subharmonic function such that $\displaystyle u(0)=0.$ So
 we can replace $\displaystyle \log \left\vert{f(z)}\right\vert
 $ by any subharmonic function $u$ in the unit disc and the "zeros
 formula" $\displaystyle \Delta \log \left\vert{f}\right\vert
 =\sum_{a\in Z(f)}{\delta _{a}}$ by the Riesz measure associated
 to $\displaystyle u,\ d\mu :=\Delta u,$ which is a positive measure.\ \par 
\quad As an application we get an alternative proof of results by Favorov
 \&  Golinskii~\cite{Golinskii12}. See also Boritchev, Golinskii
 and Kupin~\cite{BGK14}.\ \par 
\ \par 
\quad Let $\displaystyle E=\bar E\subset {\mathbb{T}}$ be a closed
 set and $\displaystyle p\geq 0,\ q>0$ real numbers ; set $\displaystyle
 \forall z\in {\mathbb{D}},\ d(z,E)$ the euclidean distance from
 $z$ to $E$ and $\varphi (z):=d(z,E)^{q}.$ Then we define a Nevanlinna
 class of functions associated to $\displaystyle E,\ p,\ q$ this
 way. For $\displaystyle p>0:$\ \par 

\begin{Dfnt}
Let $\displaystyle E=\bar E\subset {\mathbb{T}}.$ We say that
 an holomorphic function $f$ in $\displaystyle {\mathbb{D}}$
 is in the generalised Nevanlinna class $\displaystyle {\mathcal{N}}_{\varphi
 ,p}({\mathbb{D}})$ for $\displaystyle p>0$ if  $\displaystyle
 \exists \delta >0,\ \delta <1$ such that\par 
\quad \quad \quad $\displaystyle \ {\left\Vert{f}\right\Vert}_{{\mathcal{N}}_{\varphi
 ,p}}:=\sup _{1-\delta <s<1}\int_{{\mathbb{D}}}{(1-\left\vert{z}\right\vert
 )^{p-1}\varphi (sz)\log ^{+}\left\vert{f(sz)}\right\vert }<\infty .$
\end{Dfnt}
And, for $\displaystyle p=0,$\ \par 

\begin{Dfnt}
Let $\displaystyle E=\bar E\subset {\mathbb{T}}.$ We say that
 an holomorphic function $f$ is in the generalised Nevanlinna
 class $\displaystyle {\mathcal{N}}_{d(\cdot ,E)^{q},0}$ if 
 $\displaystyle \exists \delta >0,\ \delta <1$ such that\par 
\quad \quad \quad $\displaystyle \ {\left\Vert{f}\right\Vert}_{{\mathcal{N}}_{d(\cdot
 ,E)^{q},0}}:=\sup _{1-\delta <s<1}\lbrace \int_{{\mathbb{T}}}{d(se^{i\theta
 },E)^{q}\log ^{+}\left\vert{f(se^{i\theta })}\right\vert }+$\par 
\quad \quad \quad \quad \quad \quad \quad $\displaystyle \ +\int_{{\mathbb{D}}}{d(sz,E)^{q-1}\log ^{+}\left\vert{f(sz)}\right\vert
 }+\int_{{\mathbb{D}}}{(1-\left\vert{sz}\right\vert ^{2})^{q-1}\log
 ^{+}\left\vert{f(sz)}\right\vert }\rbrace <\infty .$
\end{Dfnt}
\quad And we prove the Blaschke type condition, for $\displaystyle p\geq 0,$\ \par 

\begin{Thrm}
Let $\displaystyle E=\bar E\subset {\mathbb{T}}.$ Suppose $\displaystyle
 q>0$ and $\displaystyle f\in {\mathcal{N}}_{\varphi ,p}({\mathbb{D}})$
 with $\displaystyle \ \left\vert{f(0)}\right\vert =1,$ then\par 
\quad \quad \quad $\displaystyle \ \sum_{a\in Z(f)}{(1-\left\vert{a}\right\vert
 ^{2})^{1+p}\varphi (a)}\leq c(\varphi ){\left\Vert{f}\right\Vert}_{{\mathcal{N}}_{\varphi
 ,p}}.$
\end{Thrm}
\quad As an application we get also the following results, which are
 special cases of results of Favorov \&  Golinskii~\cite{Golinskii12}.
 See also Boritchev, Golinskii and Kupin~\cite{BGK14}.\ \par 

\begin{Thrm}
Suppose that $\displaystyle f\in {\mathcal{H}}({\mathbb{D}}),\
 \left\vert{f(0)}\right\vert =1$ and\par 
\quad \quad \quad $\displaystyle \forall z\in {\mathbb{D}},\ \log ^{+}\left\vert{f(z)}\right\vert
 \leq \frac{K}{(1-\left\vert{z}\right\vert ^{2})^{p}}\frac{1}{d(z,E)^{q}},$\par 
then we have, with any $\displaystyle \epsilon >0,\ p>0,$\par 
\quad \quad \quad $\displaystyle \ \sum_{a\in Z(f)}{(1-\left\vert{a}\right\vert
 ^{2})^{1+p}d(a,E)^{(q-\alpha (E)+\epsilon )_{+}}}\leq c(p,q,\epsilon )K.$
\end{Thrm}
And in the case $\displaystyle p=0,$\ \par 

\begin{Thrm}
Suppose that $\displaystyle f\in {\mathcal{H}}({\mathbb{D}}),\
 \left\vert{f(0)}\right\vert =1$ and\par 
\quad \quad \quad $\displaystyle \forall z\in {\mathbb{D}},\ \log ^{+}\left\vert{f(z)}\right\vert
 \leq K\frac{1}{d(z,E)^{q}},$\par 
then\par 
\quad \quad \quad $\displaystyle \ \sum_{a\in Z(f)}{(1-\left\vert{a}\right\vert
 ^{2})d(a,E)^{(q-\alpha (E)+\epsilon )_{+}}}\leq c(q,\epsilon )K.$
\end{Thrm}

\subsection{Notations.}
\quad Let $\displaystyle E=\bar E\subset {\mathbb{T}}$ be a closed
 set ; we have $\displaystyle {\mathbb{T}}\backslash E=\bigcup_{j\in
 {\mathbb{N}}}{(\alpha _{j},\beta _{j})}$ where the $\displaystyle
 F_{j}:=(\alpha _{j},\beta _{j})$ are the contiguous intervals
 to $\displaystyle E.$ Set $\displaystyle 2\delta _{j}$ the length
 of the arc $\displaystyle F_{j}.$\ \par 
\quad Let $\displaystyle \Gamma _{j}:=\lbrace z=re^{i\psi }\in {\mathbb{D}}::\psi
 \in (\alpha _{j},\beta _{j})\rbrace $ the conical set based
 on $\displaystyle F_{j}$ and $\Gamma _{E}:=\lbrace z=re^{i\psi
 }\in {\mathbb{D}}::\psi \in E\rbrace .$\ \par 
Let\ \par 
\quad \quad \quad $\displaystyle \chi \in {\mathcal{C}}^{\infty }({\mathbb{R}}),\
 t\leq 2\Rightarrow \chi (t)=0,\ t\geq 3\Rightarrow \chi (t)=1.$\ \par 
\quad Now we define:\ \par 
\quad \quad \quad $\displaystyle \forall z\in \Gamma _{j},\ \psi _{j}(z):=\frac{\left\vert{z-\alpha
 _{j}}\right\vert ^{2}\left\vert{z-\beta _{j}}\right\vert ^{2}}{\delta
 _{j}^{2}},\ \eta _{j}(z):=\chi (\frac{\left\vert{z-\alpha _{j}}\right\vert
 ^{2}}{(1-\left\vert{z}\right\vert ^{2})^{2}})\chi (\frac{\left\vert{z-\beta
 _{j}}\right\vert ^{2}}{(1-\left\vert{z}\right\vert ^{2})^{2}}),$\ \par 
and\ \par 
\quad \quad \quad $\displaystyle \forall z\in \Gamma _{j},\ \varphi _{j}(z):=\eta
 _{j}(z)\psi _{j}(z)^{q}+(1-\left\vert{z}\right\vert ^{2})^{2q},\
 \forall z\in \Gamma _{E},\ \varphi _{E}(z):=(1-\left\vert{z}\right\vert
 ^{2})^{2q}.$\ \par 

\begin{Lmm}
~\label{CI25}We have\par 
\quad \quad \quad $\displaystyle \forall z\in \Gamma _{j},\ \varphi _{j}(z)\geq
 \frac{1}{3^{q}}d(z,\lbrace \alpha _{j},\beta _{j}\rbrace )^{2q}$\par 
and\par 
\quad \quad \quad $\displaystyle \forall z\in \Gamma _{j},\ \varphi _{j}(z)\leq
 (4^{q}+2^{q})d(z,\lbrace \alpha _{j},\beta _{j}\rbrace )^{2q}.$
\end{Lmm}
\quad Proof.\ \par 
We have \ \par 
\quad \quad \quad $\displaystyle \forall z\in \Gamma _{j},\ d(z,\lbrace \alpha
 _{j},\beta _{j}\rbrace )=\min (\left\vert{z-\alpha _{j}}\right\vert
 ,\left\vert{z-\beta _{j}}\right\vert ).$\ \par 
Suppose that $\displaystyle d(z,\lbrace \alpha _{j},\beta _{j}\rbrace
 )=\left\vert{z-\alpha _{j}}\right\vert $ then $\displaystyle
 \ \left\vert{z-\beta _{j}}\right\vert \geq \left\vert{z-\alpha
 _{j}}\right\vert $ hence $\displaystyle \ \left\vert{z-\beta
 _{j}}\right\vert \geq \delta _{j}.$ So\ \par 
\quad \quad \quad $\displaystyle \psi _{j}(z):=\frac{\left\vert{z-\alpha _{j}}\right\vert
 ^{2}\left\vert{z-\beta _{j}}\right\vert ^{2}}{\delta _{j}^{2}}\geq
 \left\vert{z-\alpha _{j}}\right\vert ^{2}=d(z,\lbrace \alpha
 _{j},\beta _{j}\rbrace )^{2}.$\ \par 
Now\ \par 
\quad $\bullet $ if $\eta _{j}(z)=1,$ then\ \par 
\quad \quad \quad $\displaystyle \forall z\in \Gamma _{j},\ d(z,\lbrace \alpha
 _{j},\beta _{j}\rbrace )^{2q}\leq \psi _{j}(z)^{q}\leq \eta
 _{j}(z)\psi _{j}(z)^{q}+(1-\left\vert{z}\right\vert ^{2})^{2q}=\varphi
 _{j}(z).$\ \par 
\ \par 
\quad $\bullet $ Suppose $\displaystyle d(z,\lbrace \alpha _{j},\beta
 _{j}\rbrace )=\left\vert{z-\alpha _{j}}\right\vert $ then if
 $\eta _{j}(z)<1$ then either $\displaystyle \ \left\vert{z-\alpha
 }\right\vert ^{2}\leq 3(1-\left\vert{z}\right\vert ^{2})^{2}$
 or $\displaystyle \ \left\vert{z-\beta }\right\vert ^{2}\leq
 3(1-\left\vert{z}\right\vert ^{2})^{2}.$ Suppose that $\displaystyle
 \ \left\vert{z-\alpha }\right\vert ^{2}\leq 3(1-\left\vert{z}\right\vert
 ^{2})^{2}$ we have\ \par 
\quad \quad \quad $\displaystyle d(z,\lbrace \alpha _{j},\beta _{j}\rbrace )^{2}=\left\vert{z-\alpha
 _{j}}\right\vert ^{2}\leq 3(1-\left\vert{z}\right\vert ^{2})^{2}\Rightarrow
 (1-\left\vert{z}\right\vert ^{2})^{2}\geq \frac{1}{3}d(z,\lbrace
 \alpha _{j},\beta _{j}\rbrace )^{2}$\ \par 
hence\ \par 
\quad \quad \quad $\displaystyle \varphi _{j}(z)=\eta _{j}(z)\psi _{j}(z)^{q}+(1-\left\vert{z}\right\vert
 ^{2})^{2q}\geq (1-\left\vert{z}\right\vert ^{2})^{2q}\geq \frac{1}{3^{q}}d(z,\lbrace
 \alpha _{j},\beta _{j}\rbrace )^{2q}.$\ \par 
\quad If $\displaystyle \ \left\vert{z-\beta }\right\vert ^{2}\leq
 3(1-\left\vert{z}\right\vert ^{2})^{2}$ still with $\displaystyle
 d(z,\lbrace \alpha _{j},\beta _{j}\rbrace )=\left\vert{z-\alpha
 _{j}}\right\vert $ then\ \par 
\quad \quad \quad $\displaystyle \ \left\vert{z-\alpha _{j}}\right\vert \leq \left\vert{z-\beta
 _{j}}\right\vert \leq 3(1-\left\vert{z}\right\vert ^{2})^{2}\Rightarrow
 (1-\left\vert{z}\right\vert ^{2})^{2}\geq \frac{1}{3}d(z,\lbrace
 \alpha _{j},\beta _{j}\rbrace )^{2}$\ \par 
and again $\displaystyle \varphi _{j}(z)\geq \frac{1}{3^{q}}d(z,\lbrace
 \alpha _{j},\beta _{j}\rbrace )^{2q}.$\ \par 
Hence in any cases we have $\displaystyle \varphi _{j}(z)\geq
 3^{-q}d(z,\lbrace \alpha _{j},\beta _{j}\rbrace )^{2q}.$\ \par 
\quad For the other way, still with $\displaystyle d(z,\lbrace \alpha
 _{j},\beta _{j}\rbrace )=\left\vert{z-\alpha _{j}}\right\vert
 ,$ we have, if $\eta _{j}(z)>0,$ that $\displaystyle \ \left\vert{z-\alpha
 }\right\vert ^{2}\geq 2(1-\left\vert{z}\right\vert ^{2})^{2}$
 hence, with $\displaystyle z=\rho e^{i\theta },$\ \par 
\quad \quad \quad $\displaystyle \ \left\vert{z-\beta }\right\vert ^{2}=(1-\rho
 )^{2}+\left\vert{e^{i\theta }-\beta }\right\vert ^{2}\geq \left\vert{z-\alpha
 }\right\vert ^{2}\geq 2(1-\rho ^{2})^{2}$\ \par 
hence\ \par 
\quad \quad \quad $\displaystyle \ \left\vert{e^{i\theta }-\beta }\right\vert ^{2}\geq
 (1-\rho ^{2})^{2}\Rightarrow (1-\rho )^{2}\leq \left\vert{e^{i\theta
 }-\beta }\right\vert ^{2}.$\ \par 
So\ \par 
\quad \begin{equation}  \ \left\vert{z-\beta }\right\vert ^{2}=(1-\rho
 )^{2}+\left\vert{e^{i\theta }-\beta }\right\vert ^{2}\leq \left\vert{e^{i\theta
 }-\beta }\right\vert ^{2}+\left\vert{e^{i\theta }-\beta }\right\vert
 ^{2}=2\left\vert{e^{i\theta }-\beta }\right\vert ^{2}\leq 2\delta
 ^{2}.\label{eOPC0}\end{equation}\ \par 
Putting it in $\psi $ we get\ \par 
\quad \quad \quad \begin{equation}  \psi _{j}(z):=\frac{\left\vert{z-\alpha _{j}}\right\vert
 ^{2}\left\vert{z-\beta _{j}}\right\vert ^{2}}{\delta _{j}^{2}}\leq
 2\left\vert{z-\alpha _{j}}\right\vert ^{2}=2d(z,\lbrace \alpha
 _{j},\beta _{j}\rbrace )^{2}\label{eOPC1}\end{equation}\ \par 
hence\ \par 
\quad \quad \quad $\eta _{j}(z)\psi _{j}(z)\leq 2d(z,\lbrace \alpha _{j},\beta
 _{j}\rbrace )^{2}.$\ \par 
Because $\displaystyle (1-\left\vert{z}\right\vert ^{2})^{2}\leq
 4d(z,\lbrace \alpha _{j},\beta _{j}\rbrace )^{2}$ we get\ \par 
\quad \quad \quad $\displaystyle \varphi _{j}(z)=\eta _{j}(z)\psi _{j}(z)^{q}+(1-\left\vert{z}\right\vert
 ^{2})^{2q}\leq (4^{q}+2^{q})d(z,\lbrace \alpha _{j},\beta _{j}\rbrace
 )^{2q}.$ $\hfill\blacksquare $\ \par 

\begin{Lmm}
~\label{eOPC2}There is a function $\displaystyle \varphi \in
 {\mathcal{C}}^{\infty }({\mathbb{D}})$ such that $\varphi $
 coincides with $\displaystyle \varphi _{j}$ and $\displaystyle
 \varphi _{E}$ in their domains of definition.
\end{Lmm}
\quad Proof.\ \par 
Clearly $\eta _{j}(z)\psi _{j}(z)^{q}$ is in $\displaystyle {\mathcal{C}}^{\infty
 }(\Gamma _{j})$ so the question is between $\Gamma _{j}$ and
 $\Gamma _{E}.$ But for any $\displaystyle s<1$ and $\displaystyle
 z\in \Gamma _{j}\cap D(0,s)$ we have that, for any multi index
 $\displaystyle \alpha \in {\mathbb{N}}^{2},\ \partial ^{\alpha
 }\lbrack \eta _{j}(z)\psi _{j}(z)^{q}\rbrack \rightarrow 0$
 when $\displaystyle z\rightarrow z_{0}\in \partial \Gamma _{j}\cap
 D(0,s)$ because $\displaystyle \chi (\frac{\left\vert{z-\alpha
 _{j}}\right\vert ^{2}}{(1-\left\vert{z}\right\vert ^{2})^{2}})$
 goes to $0$ with all its derivatives when $\displaystyle \ \frac{\left\vert{z-\alpha
 _{j}}\right\vert ^{2}}{(1-\left\vert{z}\right\vert ^{2})^{2}}\rightarrow
 0.$   The same for $\displaystyle \chi (\frac{\left\vert{z-\beta
 _{j}}\right\vert ^{2}}{(1-\left\vert{z}\right\vert ^{2})^{2}}).$
 So $\eta _{j}(z)\psi _{j}(z)^{q}$ extends $\displaystyle {\mathcal{C}}^{\infty
 }$ by $0$ to $\Gamma _{E}\cap D(0,s).$ And $\displaystyle \varphi
 _{E}(z):=(1-\left\vert{z}\right\vert ^{2})^{2q}$ is already
 global and $\displaystyle {\mathcal{C}}^{\infty }({\mathbb{D}}).$
 (Not in $\displaystyle {\mathcal{C}}^{\infty }(\bar {\mathbb{D}})\ !$)\ \par 
\quad So $\varphi _{j}$ being the sum of these functions extends $\displaystyle
 {\mathcal{C}}^{\infty }$ to the open disc. $\hfill\blacksquare $\ \par 
\ \par 
\quad Now we set, for $\displaystyle 0\leq s<1$ and $\displaystyle
 q>0,\ g_{s}(z):=(1-\left\vert{z}\right\vert ^{2})^{p+1}\varphi
 (sz)\in {\mathcal{C}}^{\infty }(\bar {\mathbb{D}})$ so we can
 apply the Green formula to it. Recall that $\displaystyle f_{s}(z):=f(sz).$\
 \par 
\quad In fact in the case of $\displaystyle \log \left\vert{f_{s}}\right\vert
 ,$ even if this function is not $\displaystyle {\mathcal{C}}^{2},$
 this is quite well known but for sake of completeness we give
 a proof as lemma~\ref{0_NI2}. Now, because everything works
 exactly the same way if we replace $\displaystyle \log \left\vert{f_{s}}\right\vert
 $ by $\displaystyle v(sz)$ where $v$ is a subharmonic function
 in the unit disc $\displaystyle {\mathbb{D}},$ we give also
 a proof of the Green formula in that case in lemma~\ref{0_NI3},
 in the appendix. Troughout this work we let $\displaystyle \log
 \left\vert{f}\right\vert $ instead of a general subharmonic
 function $v$ because it is the most interesting case.\ \par 
\quad With the "zero" formula: $\Delta \log \left\vert{f_{s}}\right\vert
 =\sum_{a\in Z(f_{s})}{\delta _{a}}$ we get\ \par 
\quad \quad \quad $\displaystyle \ \sum_{a\in Z(f_{s})}{g_{s}(a)}=\int_{{\mathbb{D}}}{\log
 \left\vert{f(sz)}\right\vert \triangle g_{s}(z)}+\int_{{\mathbb{T}}}{(g_{s}\partial
 _{n}\log \left\vert{f(sz)}\right\vert -\log \left\vert{f(sz)}\right\vert
 \partial _{n}g_{s})}.$\ \par 
Because $\displaystyle g_{s}=0$ on $\displaystyle {\mathbb{T}},$ we get:\ \par 
\quad \quad \quad $\displaystyle \ \sum_{a\in Z(f_{s})}{g_{s}(a)}=\int_{{\mathbb{D}}}{\log
 \left\vert{f(sz)}\right\vert \triangle g_{s}(z)}-\int_{{\mathbb{T}}}{\log
 \left\vert{f(se^{i\theta })}\right\vert \partial _{n}g_{s}(e^{i\theta
 })}.$\ \par 
If, moreover $\displaystyle p>0,\ \partial _{n}g_{s}=0$ on $\displaystyle
 {\mathbb{T}},$ hence $\displaystyle \ \sum_{a\in Z(f_{s})}{g_{s}(a)}=\int_{{\mathbb{D}}}{\log
 \left\vert{f(sz)}\right\vert \triangle g_{s}(z)}.$ So we proved\ \par 

\begin{Lmm}
~\label{eOPC3} Let $\displaystyle p>0$ and $\displaystyle f\in
 {\mathcal{H}}({\mathbb{D}})$ we have\par 
\quad \quad \quad $\displaystyle \ \sum_{a\in Z(f_{s})}{g_{s}(a)}=\int_{{\mathbb{D}}}{\log
 \left\vert{f(sz)}\right\vert \triangle g_{s}(z)}.$
\end{Lmm}
\quad We have to compute\ \par 
\quad \quad \quad $\displaystyle \triangle g_{s}(z)\log \left\vert{f(sz)}\right\vert
 =\triangle g_{s}(z)\log ^{+}\left\vert{f(sz)}\right\vert -\triangle
 g_{s}(z)\log ^{-}\left\vert{f(sz)}\right\vert .$\ \par 
We have $\Delta g_{s}=4\bar \partial \partial g_{s}$ hence\ \par 
\quad \quad \quad $\displaystyle \Delta g_{s}(z)=\Delta \lbrack (1-\left\vert{z}\right\vert
 ^{2})^{p+1}\varphi (sz)\rbrack =\varphi (sz)\Delta \lbrack (1-\left\vert{z}\right\vert
 ^{2})^{p+1}\rbrack +(1-\left\vert{z}\right\vert ^{2})^{p+1}\Delta
 \lbrack \varphi (sz)\rbrack +$\ \par 
\quad \quad \quad \quad \quad \quad \quad $\displaystyle +8\Re \lbrack \partial ((1-\left\vert{z}\right\vert
 ^{2})^{p+1})\bar \partial (\varphi (sz))\rbrack .$\ \par 
\ \par 
\quad Recall that, with $\displaystyle \varphi _{A,j}(z):=\eta _{j}(z)\psi
 _{j}(z)^{q}$ and $\displaystyle \varphi _{C,j}:=(1-\left\vert{z}\right\vert
 ^{2})^{2q},$\ \par 
\quad \quad \quad $\displaystyle \forall z\in \Gamma _{j},\ \varphi _{j}(z):=\varphi
 _{A,j}(z)+\varphi _{C,j}(z)\ ;$\ \par 
we start with the last term.\ \par 

\section{Estimates on $\varphi _{C,j}(z):=(1-\left\vert{z}\right\vert
 ^{2})^{2q}.$}
\quad In this case\ \par 
\quad \quad \quad $\displaystyle g_{C,s}(z):=(1-\left\vert{z}\right\vert ^{2})^{p+1}\varphi
 _{C}(sz)=(1-\left\vert{z}\right\vert ^{2})^{p+1}(1-\left\vert{sz}\right\vert
 ^{2})^{2q}.$\ \par 
So we have to compute $\Delta \lbrack (1-\left\vert{z}\right\vert
 ^{2})^{p+1}\varphi _{C}(z)\rbrack =A_{1}+A_{2}+A_{3}$ with:\ \par 
\quad \quad \quad $\displaystyle A_{1}:=(1-\left\vert{sz}\right\vert ^{2})^{2q}\Delta
 ((1-\left\vert{z}\right\vert ^{2})^{p+1})=$\ \par 
\quad \quad \quad \quad \quad $\displaystyle =-4(p+1)(1-\left\vert{z}\right\vert ^{2})^{p}(1-\left\vert{sz}\right\vert
 ^{2})^{2q}+4p(p+1)(1-\left\vert{z}\right\vert ^{2})^{p-1}\left\vert{z}\right\vert
 ^{2}(1-\left\vert{sz}\right\vert ^{2})^{2q}\ ;$\ \par 
\quad \quad \quad $\displaystyle A_{2}:=(1-\left\vert{z}\right\vert ^{2})^{p+1}\Delta
 ((1-\left\vert{sz}\right\vert ^{2})^{2q})=$\ \par 
\quad \quad \quad \quad \quad $\displaystyle -8sq(1-\left\vert{z}\right\vert ^{2})^{p+1}(1-\left\vert{sz}\right\vert
 ^{2})^{2q-1}+8q(2q-1)(1-\left\vert{z}\right\vert ^{2})^{p+1}\left\vert{z}\right\vert
 ^{2}(1-\left\vert{sz}\right\vert ^{2})^{2q-2}\ ;$\ \par 
\quad \quad \quad $\displaystyle A_{3}:=8\Re \lbrack \partial ((1-\left\vert{z}\right\vert
 ^{2})^{p+1})\bar \partial ((1-\left\vert{sz}\right\vert ^{2})^{2q})\rbrack
 =$\ \par 
\quad \quad \quad \quad \quad $\displaystyle =16sq(p+1)(1-\left\vert{z}\right\vert ^{2})^{p}(1-\left\vert{sz}\right\vert
 ^{2})^{2q-1}\left\vert{z}\right\vert ^{2}.$\ \par 
\quad We shall consider the terms  $\displaystyle \Delta g_{C,s}(z)\log
 ^{+}\left\vert{f(sz)}\right\vert .$ We shall use\ \par 

\begin{Lmm}
~\label{CI0} For $\displaystyle p>0$ we have:\par 
\quad \quad \quad $\displaystyle \forall z\in {\mathbb{D}},\ \Delta g_{C,s}(z)\leq
 c(p,q)(1-\left\vert{z}\right\vert ^{2})^{p-1}\left\vert{z}\right\vert
 ^{2}(1-\left\vert{sz}\right\vert ^{2})^{2q}.$\par 
And for $\displaystyle p=0$ we have:\par 
\quad \quad \quad $\displaystyle \forall z\in {\mathbb{D}},\ \Delta g_{C,s}(z)\leq
 c(q)\left\vert{z}\right\vert ^{2}(1-\left\vert{sz}\right\vert
 ^{2})^{2q-1},$\par 
with $\displaystyle c(q):=8q(2q-1)+16q(p+1)$ (hence $\displaystyle c(0)=0$).
\end{Lmm}
\quad Proof.\ \par 
We have\ \par 
\quad \quad \quad $\displaystyle A_{1}\leq 4p(p+1)(1-\left\vert{z}\right\vert ^{2})^{p-1}\left\vert{z}\right\vert
 ^{2}(1-\left\vert{sz}\right\vert ^{2})^{2q},$\ \par 
because $\displaystyle -4(p+1)(1-\left\vert{z}\right\vert ^{2})^{p}(1-\left\vert{sz}\right\vert
 ^{2})^{2q}$ is negative.\ \par 
\quad \quad \quad $\displaystyle A_{2}\leq 8q(2q-1)(1-\left\vert{z}\right\vert
 ^{2})^{p+1}\left\vert{z}\right\vert ^{2}(1-\left\vert{sz}\right\vert
 ^{2})^{2q-2},$\ \par 
because $\displaystyle -8sq(1-\left\vert{z}\right\vert ^{2})^{p+1}(1-\left\vert{sz}\right\vert
 ^{2})^{2q-1}$ is negative. So adding, we get\ \par 
\quad \quad \quad $\displaystyle \Delta g_{C,s}(z)\leq 4p(p+1)(1-\left\vert{z}\right\vert
 ^{2})^{p-1}\left\vert{z}\right\vert ^{2}(1-\left\vert{sz}\right\vert
 ^{2})^{2q}+$\ \par 
\quad \quad \quad \quad \quad \quad \quad \quad \quad \quad \quad $\displaystyle +8q(2q-1)(1-\left\vert{z}\right\vert ^{2})^{p+1}\left\vert{z}\right\vert
 ^{2}(1-\left\vert{sz}\right\vert ^{2})^{2q-2}+$\ \par 
\quad \quad \quad \quad \quad \quad \quad \quad \quad \quad \quad $\displaystyle +16sq(p+1)(1-\left\vert{z}\right\vert ^{2})^{p}(1-\left\vert{sz}\right\vert
 ^{2})^{2q-1}\left\vert{z}\right\vert ^{2}.$\ \par 
If $\displaystyle p>0$ we use $\displaystyle (1-\left\vert{z}\right\vert
 ^{2})\leq (1-\left\vert{sz}\right\vert ^{2})$ to get\ \par 
\quad \quad \quad $\displaystyle (1-\left\vert{z}\right\vert ^{2})^{p+1}\left\vert{z}\right\vert
 ^{2}(1-\left\vert{sz}\right\vert ^{2})^{2q-2}\leq (1-\left\vert{z}\right\vert
 ^{2})^{p-1}\left\vert{z}\right\vert ^{2}(1-\left\vert{sz}\right\vert
 ^{2})^{2q},$\ \par 
and\ \par 
\quad \quad \quad $\displaystyle (1-\left\vert{z}\right\vert ^{2})^{p}(1-\left\vert{sz}\right\vert
 ^{2})^{2q-1}\left\vert{z}\right\vert ^{2}\leq (1-\left\vert{z}\right\vert
 ^{2})^{p-1}\left\vert{z}\right\vert ^{2}(1-\left\vert{sz}\right\vert
 ^{2})^{2q}.$\ \par 
If $\displaystyle p=0$ we keep\ \par 
\quad \quad \quad $\displaystyle (1-\left\vert{z}\right\vert ^{2})\left\vert{z}\right\vert
 ^{2}(1-\left\vert{sz}\right\vert ^{2})^{2q-2}\leq \left\vert{z}\right\vert
 ^{2}(1-\left\vert{sz}\right\vert ^{2})^{2q-1},$\ \par 
So, setting, for $\displaystyle p>0,$\ \par 
\quad \quad \quad $\displaystyle c(p,q):=4p(p+1)+8q(2q-1)+16q(p+1),$\ \par 
and\ \par 
\quad \quad \quad $\displaystyle c(q):=8q(2q-1)+16q(p+1),$\ \par 
which ends the proof of the lemma. $\hfill\blacksquare $\ \par 

\begin{Prps}
~\label{CI7}We have, for $\displaystyle p>0,$\par 
\quad \quad \quad $\displaystyle \ \int_{{\mathbb{D}}}{\Delta g_{C,s}(z)\log ^{+}\left\vert{f(sz)}\right\vert
 }\leq c(p,q)\int_{{\mathbb{D}}}{(1-\left\vert{z}\right\vert
 ^{2})^{p-1}\varphi _{C}(sz)\log ^{+}\left\vert{f(sz)}\right\vert }.$\par 
And for $\displaystyle p=0,$\par 
\quad \quad \quad $\displaystyle \ \int_{{\mathbb{D}}}{\Delta g_{C,s}(z)\log ^{+}\left\vert{f(sz)}\right\vert
 }\leq c(q)\int_{{\mathbb{D}}}{(1-\left\vert{sz}\right\vert ^{2})^{2q-1}\log
 ^{+}\left\vert{f(sz)}\right\vert }$\par 
with $\displaystyle c(0)=0.$
\end{Prps}
\quad Proof.\ \par 
Integrating the estimates of lemma~\ref{CI0} we get the proposition.
 $\hfill\blacksquare $\ \par 
\quad In order to consider the terms containing $\displaystyle \log
 ^{-}\left\vert{f(sz)}\right\vert $ we shall need:\ \par 

\begin{Lmm}
~\label{CI3} We have, for $\displaystyle p\geq 0$ and any $\displaystyle
 s\geq 1/2,$\par 
\quad \quad \quad $\displaystyle \forall z\in {\mathbb{D}},\ -\Delta g_{C,s}(z)\leq
 -4p(p+1)(1-\left\vert{z}\right\vert ^{2})^{p-1}\left\vert{z}\right\vert
 ^{2}(1-\left\vert{sz}\right\vert ^{2})^{2q}+$\par 
\quad \quad \quad \quad \quad \quad \quad \quad \quad \quad \quad \quad \quad \quad \quad \quad \quad $\displaystyle +4\lbrack (p+1)+2sq\rbrack (1-\left\vert{z}\right\vert
 ^{2})^{p}(1-\left\vert{sz}\right\vert ^{2})^{2q}.$
\end{Lmm}
\quad Proof.\ \par 
With $\displaystyle -\Delta g_{C,s}(z)$ we get:\ \par 
\quad \quad \quad $\displaystyle -A_{1}\leq 4(p+1)(1-\left\vert{z}\right\vert ^{2})^{p}(1-\left\vert{sz}\right\vert
 ^{2})^{2q}-$\ \par 
\quad \quad \quad \quad \quad $\displaystyle -4p(p+1)(1-\left\vert{z}\right\vert ^{2})^{p-1}\left\vert{z}\right\vert
 ^{2}(1-\left\vert{sz}\right\vert ^{2})^{2q}\ ;$\ \par 
\quad \quad \quad $\displaystyle -A_{2}=8sq(1-\left\vert{z}\right\vert ^{2})^{p+1}(1-\left\vert{sz}\right\vert
 ^{2})^{2q-1}-$\ \par 
\quad \quad \quad \quad \quad $\displaystyle -8q(2q-1)(1-\left\vert{z}\right\vert ^{2})^{p+1}\left\vert{z}\right\vert
 ^{2}(1-\left\vert{sz}\right\vert ^{2})^{2q-2}\ ;$\ \par 
and,\ \par 
\quad \quad \quad $\displaystyle -A_{3}=-16sq(p+1)(1-\left\vert{z}\right\vert ^{2})^{p}(1-\left\vert{sz}\right\vert
 ^{2})^{2q-1}\left\vert{z}\right\vert ^{2}.$\ \par 
We have two cases:\ \par 
\quad $\bullet $ $\displaystyle 2q-1\geq 0$ then $\displaystyle -A_{2}\leq
 8sq(1-\left\vert{z}\right\vert ^{2})^{p+1}(1-\left\vert{sz}\right\vert
 ^{2})^{2q-1}.$\ \par 
\quad $\bullet $ $\displaystyle 2q-1<0$ then:\ \par 
\quad \quad \quad $\displaystyle 8q(1-2q)(1-\left\vert{z}\right\vert ^{2})^{p}\left\vert{z}\right\vert
 ^{2}(1-\left\vert{sz}\right\vert ^{2})^{2q-1}-$\ \par 
\quad \quad \quad $\displaystyle -16sq(p+1)(1-\left\vert{z}\right\vert ^{2})^{p}(1-\left\vert{sz}\right\vert
 ^{2})^{2q-1}\left\vert{z}\right\vert ^{2}=$\ \par 
\quad \quad \quad \quad \quad \quad $\displaystyle =-\lbrack 16sq(p+1)-8q(1-2q)\rbrack (1-\left\vert{z}\right\vert
 ^{2})^{p}(1-\left\vert{sz}\right\vert ^{2})^{2q-1}\left\vert{z}\right\vert
 ^{2}.$\ \par 
But\ \par 
\quad \quad \quad $\displaystyle \lbrack 16sq(p+1)-8q(1-2q)\rbrack =16sq-8q+16sqp+16q^{2}\geq
 8q(2s-1)\geq 0$\ \par 
provided that $\displaystyle s\geq 1/2.$\ \par 
So in any cases, with $\displaystyle s\geq 1/2,$ we get for $\displaystyle
 p>0,$\ \par 
\quad \quad \quad $\displaystyle \forall z\in {\mathbb{D}},\ -\Delta g_{C,s}(z)\leq
 -4p(p+1)(1-\left\vert{z}\right\vert ^{2})^{p-1}\left\vert{z}\right\vert
 ^{2}(1-\left\vert{sz}\right\vert ^{2})^{2q}+$\ \par 
\quad \quad \quad \quad \quad \quad \quad \quad \quad \quad \quad \quad \quad \quad \quad \quad \quad $\displaystyle +4(p+1)(1-\left\vert{z}\right\vert ^{2})^{p}(1-\left\vert{sz}\right\vert
 ^{2})^{2q}+$\ \par 
\quad \quad \quad \quad \quad \quad \quad \quad \quad \quad \quad \quad \quad \quad \quad \quad \quad $\displaystyle +8sq(1-\left\vert{z}\right\vert ^{2})^{p+1}(1-\left\vert{sz}\right\vert
 ^{2})^{2q-1}\leq $\ \par 
\quad \quad \quad \quad \quad \quad \quad \quad \quad \quad \quad \quad $\displaystyle \leq -4p(p+1)(1-\left\vert{z}\right\vert ^{2})^{p-1}\left\vert{z}\right\vert
 ^{2}(1-\left\vert{sz}\right\vert ^{2})^{2q}+$\ \par 
\quad \quad \quad \quad \quad \quad \quad \quad \quad \quad \quad \quad \quad \quad \quad \quad \quad $\displaystyle +4\lbrack (p+1)+2sq\rbrack (1-\left\vert{z}\right\vert
 ^{2})^{p}(1-\left\vert{sz}\right\vert ^{2})^{2q}.$\ \par 
And for $\displaystyle p=0,$\ \par 
\quad \quad \quad $\displaystyle \forall z\in {\mathbb{D}},\ -\Delta g_{C,s}(z)\leq
 4\lbrack 1+2sq\rbrack (1-\left\vert{sz}\right\vert ^{2})^{2q}.$\ \par 
So we proved the lemma. $\hfill\blacksquare $\ \par 

\begin{Prps}
~\label{CI4}We have with $\displaystyle \ \left\vert{f(0)}\right\vert
 =1$ and $\displaystyle p\geq 0,$\par 
\quad \quad \quad $\displaystyle -\int_{{\mathbb{D}}}{\Delta g_{C,s}(z)\log ^{-}\left\vert{f(sz)}\right\vert
 }\leq 4\lbrack (p+1)+2sq\rbrack \int_{{\mathbb{D}}}{(1-\left\vert{z}\right\vert
 ^{2})^{p}(1-\left\vert{sz}\right\vert ^{2})^{2q}\log ^{+}\left\vert{f(sz)}\right\vert
 }.$
\end{Prps}
\quad Proof.\ \par 
Passing in polar coordinates we get\ \par 
\quad $\displaystyle \ \int_{{\mathbb{D}}}{(1-\left\vert{z}\right\vert
 ^{2})^{p}(1-\left\vert{sz}\right\vert ^{2})^{2q}\log ^{-}\left\vert{f(sz)}\right\vert
 }=\int_{0}^{1}{(1-\rho ^{2})^{p}(1-s^{2}\rho ^{2})^{2q}\lbrace
 \int_{{\mathbb{T}}}{\log ^{-}\left\vert{f(s\rho e^{i\theta })}\right\vert
 }\rbrace \rho d\rho }\ ;$\ \par 
by the subharmonicity of $\displaystyle \log \left\vert{f(sz)}\right\vert
 $ and the fact $\displaystyle \ \left\vert{f(0)}\right\vert =1,$ we get\ \par 
\quad \quad \quad $\displaystyle \ \int_{{\mathbb{T}}}{\log ^{-}\left\vert{f(s\rho
 e^{i\theta })}\right\vert }\leq \int_{{\mathbb{T}}}{\log ^{+}\left\vert{f(s\rho
 e^{i\theta })}\right\vert }$\ \par 
so\ \par 
\quad \quad \quad $\displaystyle \ \int_{{\mathbb{D}}}{(1-\left\vert{z}\right\vert
 ^{2})^{p}(1-\left\vert{sz}\right\vert ^{2})^{2q}\log ^{-}\left\vert{f(sz)}\right\vert
 }\leq \int_{{\mathbb{D}}}{(1-\left\vert{z}\right\vert ^{2})^{p}(1-\left\vert{sz}\right\vert
 ^{2})^{2q}\log ^{+}\left\vert{f(sz)}\right\vert }.$\ \par 
Now using lemma~\ref{CI3},\ \par 
\quad \quad \quad $\displaystyle \ -\int_{{\mathbb{D}}}{\Delta g_{C,s}(z)\log ^{-}\left\vert{f(sz)}\right\vert
 }\leq -4p(p+1)\int_{{\mathbb{D}}}{(1-\left\vert{z}\right\vert
 ^{2})^{p-1}\left\vert{z}\right\vert ^{2}(1-\left\vert{sz}\right\vert
 ^{2})^{2q}\log ^{-}\left\vert{f(sz)}\right\vert }+$\ \par 
\quad \quad \quad \quad \quad \quad \quad \quad \quad \quad \quad \quad \quad \quad \quad \quad \quad $\displaystyle +4\lbrack (p+1)+2sq\rbrack \int_{{\mathbb{D}}}{(1-\left\vert{z}\right\vert
 ^{2})^{p}(1-\left\vert{sz}\right\vert ^{2})^{2q}\log ^{-}\left\vert{f(sz)}\right\vert
 }\leq $\ \par 
\quad \quad \quad \quad \quad \quad \quad \quad \quad \quad \quad \quad \quad \quad \quad $\displaystyle \leq 4\lbrack (p+1)+2sq\rbrack \int_{{\mathbb{D}}}{(1-\left\vert{z}\right\vert
 ^{2})^{p}(1-\left\vert{sz}\right\vert ^{2})^{2q}\log ^{-}\left\vert{f(sz)}\right\vert
 }\leq $\ \par 
\quad \quad \quad \quad \quad \quad \quad \quad \quad \quad \quad \quad \quad \quad \quad $\displaystyle \leq 4\lbrack (p+1)+2sq\rbrack \int_{{\mathbb{D}}}{(1-\left\vert{z}\right\vert
 ^{2})^{p}(1-\left\vert{sz}\right\vert ^{2})^{2q}\log ^{+}\left\vert{f(sz)}\right\vert
 }.$\ \par 
This ends the proof of the proposition. $\hfill\blacksquare $\ \par 

\section{Estimates on $\varphi _{A,j}(z):=\psi _{j}(sz)^{q}\eta _{j}(sz).$}
\quad We set\ \par 
\quad \quad \quad $\displaystyle g_{A,s}(z):=(1-\left\vert{z}\right\vert ^{2})^{p+1}\sum_{j\in
 {\mathbb{N}}}{{\11}_{\Gamma _{j}}(z)\varphi _{A,j}(sz)},$\ \par 
and we have seen that $\displaystyle g_{A,s}(z)\in {\mathcal{C}}^{\infty
 }({\mathbb{D}}.)$\ \par 
We shall compute\ \par 
\quad \quad \quad $\displaystyle \triangle g_{A,s}(z)\log \left\vert{f(sz)}\right\vert
 =\triangle g_{A,s}(z)\log ^{+}\left\vert{f(sz)}\right\vert -\triangle
 g_{A,s}(z)\log ^{-}\left\vert{f(sz)}\right\vert .$\ \par 
We have $\Delta g_{A,s}=4\bar \partial \partial g_{A,s}$ hence here:\ \par 
\quad \quad \quad $\displaystyle \forall z\in \Gamma _{j},\ \Delta \lbrack (1-\left\vert{z}\right\vert
 ^{2})^{p+1}\varphi _{A,j}(sz)\rbrack =\varphi _{A,j}(sz)\Delta
 \lbrack (1-\left\vert{z}\right\vert ^{2})^{p+1}\rbrack +$\ \par 
\quad \quad \quad \quad \quad \quad \quad \quad \quad \quad \quad \quad \quad \quad \quad \quad \quad \quad \quad \quad \quad $\displaystyle +(1-\left\vert{z}\right\vert ^{2})^{p+1}\Delta
 \lbrack \varphi _{A,j}(sz)\rbrack +$\ \par 
\quad \quad \quad \quad \quad \quad \quad \quad \quad \quad \quad \quad \quad \quad \quad \quad \quad \quad \quad \quad \quad $\displaystyle +8\Re \lbrack \partial ((1-\left\vert{z}\right\vert
 ^{2})^{p+1})\bar \partial (\varphi _{A,j}(sz))\rbrack =:$\ \par 
\quad \quad \quad \quad \quad \quad \quad \quad \quad \quad \quad \quad \quad \quad \quad \quad \quad \quad $\displaystyle =:A_{1}+A_{2}+A_{3}.$\ \par 

\begin{Lmm}
~\label{eOP2}We have:\par 
\quad \quad \quad $\displaystyle A_{1}:=\varphi _{A,j}(sz)\Delta \lbrack (1-\left\vert{z}\right\vert
 ^{2})^{p+1}\rbrack =4p(p+1)(1-\left\vert{z}\right\vert ^{2})^{p-1}\left\vert{z}\right\vert
 ^{2}\varphi _{A,j}(sz)-$\par 
\quad \quad \quad \quad \quad \quad \quad $\displaystyle -4(p+1)(1-\left\vert{z}\right\vert ^{2})^{p}\varphi
 _{A,j}(sz)=:$\par 
\quad \quad \quad \quad \quad $\displaystyle =:A'_{1}-A''_{1}$\par 
with, for $\displaystyle z\in \Gamma _{j},$\par 
\quad \quad \quad $\displaystyle A'_{1}:=4p(p+1)(1-\left\vert{z}\right\vert ^{2})^{p-1}\left\vert{z}\right\vert
 ^{2}\varphi _{A,j}(sz)\leq 2^{2q}{\times}4p(p+1)(1-\left\vert{z}\right\vert
 ^{2})^{p-1}\left\vert{z}\right\vert ^{2}d(sz,E)^{2q}.$\par 
\quad \quad \quad $\displaystyle A"_{1}:=4(p+1)(1-\left\vert{z}\right\vert ^{2})^{p}\varphi
 _{A,j}(sz)\leq 2^{2q}{\times}4(p+1)(1-\left\vert{z}\right\vert
 ^{2})^{p}d(sz,E)^{2q}.$
\end{Lmm}
\quad Proof.\ \par 
A simple computation of $\displaystyle \Delta \lbrack (1-\left\vert{z}\right\vert
 ^{2})^{p+1}\rbrack =4\partial \bar \partial \lbrack (1-\left\vert{z}\right\vert
 ^{2})^{p+1}\rbrack $ with lemma~\ref{CI25} gives the result.
 $\hfill\blacksquare $\ \par 
\ \par 
\quad We have, just using $\Delta =4\partial \bar \partial ,$\ \par 
\quad \quad \quad $\displaystyle A_{2}:=(1-\left\vert{z}\right\vert ^{2})^{p+1}\Delta
 \lbrack \varphi _{A,j}(sz)\rbrack =(1-\left\vert{z}\right\vert
 ^{2})^{p+1}\Delta \lbrack \eta _{j}(sz)\psi _{j}(sz)^{q}\rbrack =$\ \par 
\quad \quad \quad \quad \quad \quad \quad $\displaystyle =(1-\left\vert{z}\right\vert ^{2})^{p+1}\eta _{j}(sz)\Delta
 \lbrack \psi _{j}(sz)^{q}\rbrack +$\ \par 
\quad \quad \quad \quad \quad \quad \quad \quad \quad \quad \quad $\displaystyle +(1-\left\vert{z}\right\vert ^{2})^{p+1}\psi _{j}(sz)^{q}\Delta
 \lbrack \eta _{j}(sz)\rbrack +$\ \par 
\quad \quad \quad \quad \quad \quad \quad \quad \quad \quad \quad $\displaystyle +(1-\left\vert{z}\right\vert ^{2})^{p+1}8\Re \lbrack
 \partial (\eta _{j}(sz))\bar \partial (\psi _{j}(sz)^{q})\rbrack
 =:A_{2,1}+A_{2,2}+A_{2,3}.$\ \par 

\begin{Lmm}
~\label{eOP3}We have\par 
$\displaystyle \forall z\in \Gamma _{j},\ A_{2,1}(s,z):=(1-\left\vert{z}\right\vert
 ^{2})^{p+1}\eta _{j}(sz)\Delta \lbrack \psi _{j}(sz)^{q}\rbrack =$\par 
$\displaystyle =4q^{2}(1-\left\vert{z}\right\vert ^{2})^{p+1}\eta
 _{j}(sz)\frac{\left\vert{sz-\alpha _{j}}\right\vert ^{2q-2}\left\vert{sz-\beta
 _{j}}\right\vert ^{2q-2}}{\delta _{j}^{2q}}\lbrace \left\vert{sz-\beta
 _{j}}\right\vert ^{2}+\left\vert{sz-\alpha _{j}}\right\vert
 ^{2}+2\Re \lbrack (sz-\alpha _{j})(\bar z-\bar \beta _{j})\rbrack
 \rbrace .$\par 
Hence $\displaystyle \forall z\in \Gamma _{j},\ A_{2,1}(s,z)\geq
 0$ and, for $\displaystyle s\geq 1/2,$\par 
\quad \quad \quad $\displaystyle \forall z\in \Gamma _{j},\ A_{2,1}(s,z)\leq 4^{2q+2}q^{2}(1-\left\vert{z}\right\vert
 ^{2})^{p}cd(sz,E)^{2q-1}.$
\end{Lmm}
\quad Proof.\ \par 
We just apply lemma~\ref{CI23} with $\Delta =4\partial \bar \partial
 ,$ to get the first assertion. Then we apply remark~\ref{CI24}
 to get $\displaystyle \forall z\in \Gamma _{j},\ A_{2,1}(z)\geq
 0.$ Now for the third assertion we notice that $\eta _{j}(sz)\leq
 1$ then, using~(\ref{eOPC0}) in lemma~\ref{CI25} with $0<\eta
 _{j}(z),$  we get $\displaystyle \ \left\vert{sz-\alpha _{j}}\right\vert
 \leq {\sqrt{2}}\delta _{j}$ and $\displaystyle \ \left\vert{sz-\beta
 _{j}}\right\vert \leq {\sqrt{2}}\delta _{j},$ hence\ \par 
\quad $\displaystyle \forall z\in \Gamma _{j},\ \left\vert{A_{2,1}(s,z)}\right\vert
 \leq 8{\times}4^{2q}q^{2}(1-\left\vert{z}\right\vert ^{2})^{p+1}cd(sz,E)^{2q}\lbrace
 \left\vert{sz-\beta _{j}}\right\vert ^{-2}+\left\vert{sz-\alpha
 _{j}}\right\vert ^{-2}\rbrace $\ \par 
using lemma~\ref{CI25}. But $\displaystyle (1-\left\vert{z}\right\vert
 ^{2})\leq 2\left\vert{z-\gamma }\right\vert $ for any $\displaystyle
 \gamma \in {\mathbb{T}}$ so we get\ \par 
\quad \quad \quad $\displaystyle \forall z\in \Gamma _{j},\ \left\vert{A_{2,1}(s,z)}\right\vert
 \leq 4^{2q+2}q^{2}(1-\left\vert{z}\right\vert ^{2})^{p}cd(sz,E)^{2q-1},$\ \par 
which ends the proof of the lemma. $\hfill\blacksquare $\ \par 
\ \par 
\quad We set $\displaystyle \chi _{\alpha }(z):=\chi (\frac{\left\vert{z-\alpha
 }\right\vert ^{2}}{(1-\left\vert{z}\right\vert ^{2})^{2}}),\
 \chi _{\beta }(z):=\chi (\frac{\left\vert{z-\beta }\right\vert
 ^{2}}{(1-\left\vert{z}\right\vert ^{2})^{2}})$ and we set $\displaystyle
 \ \left\vert{\chi '}\right\vert :=\max (\left\vert{\chi _{\alpha
 }'}\right\vert ,\left\vert{\chi _{\beta }'}\right\vert )$ and
 $\displaystyle \ \left\vert{\chi ''}\right\vert :=\max (\left\vert{\chi
 _{\alpha }''}\right\vert ,\left\vert{\chi _{\beta }''}\right\vert ).$\ \par 

\begin{Lmm}
~\label{CI26}We have\par 
\quad \quad \quad $\displaystyle \forall z\in \Gamma _{j},\ A_{2,2}(s,z):=(1-\left\vert{z}\right\vert
 ^{2})^{p+1}\psi _{j}(sz)^{q}\Delta \lbrack \eta _{j}(sz)\rbrack
 \Rightarrow $\par 
\quad \quad \quad \quad \quad $\displaystyle \Rightarrow \left\vert{A_{2,2}(s,z)}\right\vert
 \lesssim (\left\vert{\chi '}\right\vert +\left\vert{\chi ''}\right\vert
 )(1-\left\vert{z}\right\vert ^{2})^{p}\psi _{j}(sz)^{q-1/2}\lesssim $\par 
\quad \quad \quad \quad \quad \quad \quad \quad \quad \quad \quad $\displaystyle \lesssim (\left\vert{\chi '}\right\vert +\left\vert{\chi
 ''}\right\vert )(1-\left\vert{z}\right\vert ^{2})^{p}d(sz,E)^{2q-1}.$\par 
And, for $\displaystyle p>0,$\par 
\quad \quad \quad $\displaystyle \ \left\vert{A_{2,2}(s,z)}\right\vert \lesssim
 (\left\vert{\chi '}\right\vert +\left\vert{\chi ''}\right\vert
 )(1-\left\vert{z}\right\vert ^{2})^{p-1}d(sz,E)^{2q}.$
\end{Lmm}
\quad Proof.\ \par 
We have \ \par 
\quad \quad \quad $\displaystyle \partial \bar \partial \lbrack \eta _{j}(sz)\rbrack
 =\chi _{\alpha }(z)\partial \bar \partial \lbrack \chi _{\beta
 }(z)\rbrack +\chi _{\beta }(z)\partial \bar \partial \lbrack
 \chi _{\alpha }(z)\rbrack +2\Re \lbrack \partial \chi _{\alpha
 }(z)\bar \partial \lbrack \chi _{\beta }(z)\rbrack .$\ \par 
and by lemma~\ref{eOP0}:\ \par 
\quad \quad \quad $\displaystyle \ \left\vert{\bar \partial \lbrack \chi _{\alpha
 }(z)\rbrack }\right\vert \leq 3\left\vert{\chi '()}\right\vert
 (\lambda +1)(1-\left\vert{z}\right\vert ^{2})^{-1}.$\ \par 
\quad \quad \quad $\displaystyle \ \left\vert{\partial \bar \partial \chi _{\beta
 }}\right\vert \lesssim (\left\vert{\chi '}\right\vert +\left\vert{\chi
 ''}\right\vert )(1-\left\vert{z}\right\vert ^{2})^{-2}.$\ \par 
So\ \par 
\quad \quad \quad $\displaystyle \ \left\vert{\Delta \eta _{j}}\right\vert \lesssim
 (\left\vert{\chi '}\right\vert +\left\vert{\chi ''}\right\vert
 )(1-\left\vert{z}\right\vert ^{2})^{-2}$\ \par 
hence\ \par 
\quad \quad \quad $\displaystyle \ \left\vert{A_{2,2}(s,z)}\right\vert =(1-\left\vert{z}\right\vert
 ^{2})^{p+1}\psi _{j}(sz)^{q}\left\vert{\Delta \lbrack \eta _{j}(sz)\rbrack
 }\right\vert \lesssim $\ \par 
\quad \quad \quad \quad \quad \quad \quad \quad \quad $\displaystyle \lesssim (1-\left\vert{z}\right\vert ^{2})^{p+1}\psi
 _{j}(sz)^{q}(\left\vert{\chi '}\right\vert +\left\vert{\chi
 ''}\right\vert )(1-\left\vert{sz}\right\vert ^{2})^{-2}.$\ \par 
Because $\displaystyle (1-\left\vert{z}\right\vert ^{2})\leq
 (1-\left\vert{sz}\right\vert ^{2})$ we get\ \par 
\quad \quad \quad $\displaystyle \ \left\vert{A_{2,2}(s,z)}\right\vert \lesssim
 (\left\vert{\chi '}\right\vert +\left\vert{\chi ''}\right\vert
 )(1-\left\vert{z}\right\vert ^{2})^{p-1}\psi _{j}(sz)^{q}.$\ \par 
\quad Now by lemma~\ref{eOP1} we get, if $\Delta \eta _{j}\neq 0,$\ \par 
\quad \quad \quad $\displaystyle \forall z\in \Gamma _{j},\ 2(1-\left\vert{z}\right\vert
 ^{2})^{2}\leq \psi _{j}(z)\leq 3(1-\left\vert{z}\right\vert ^{2})^{2}$\ \par 
and\ \par 
\quad \quad \quad $\displaystyle \ \left\vert{A_{2,2}(s,z)}\right\vert \lesssim
 (\left\vert{\chi '}\right\vert +\left\vert{\chi ''}\right\vert
 )(1-\left\vert{z}\right\vert ^{2})^{p}\psi _{j}(sz)^{q-1/2}.$\ \par 
Because $\displaystyle (1-\left\vert{z}\right\vert ^{2})\leq
 d(z,E)$ we get, for $\displaystyle p>0,\ (1-\left\vert{z}\right\vert
 ^{2})^{p}\psi _{j}(sz)^{q-1/2}\leq (1-\left\vert{z}\right\vert
 ^{2})^{p-1}d(sz,E)^{2q}.$\ \par 
It remains to use lemma~\ref{CI25} to get the result. $\hfill\blacksquare
 $\ \par 

\begin{Lmm}
We have\par 
\quad \quad \quad $\displaystyle \forall z\in \Gamma _{j},\ A_{2,3}:=(1-\left\vert{z}\right\vert
 ^{2})^{p+1}8\Re \lbrack \partial (\eta _{j}(sz))\bar \partial
 (\psi _{j}(sz)^{q})\rbrack \Rightarrow $\par 
\quad \quad \quad \quad \quad \quad \quad $\displaystyle \Rightarrow \left\vert{A_{2,3}}\right\vert \lesssim
 \left\vert{\chi '}\right\vert (1-\left\vert{z}\right\vert ^{2})^{p}\psi
 _{j}(sz)^{q-1/2}\lesssim \left\vert{\chi '}\right\vert (1-\left\vert{z}\right\vert
 ^{2})^{p}d(sz,E)^{2q-1}.$
\end{Lmm}
\quad Proof.\ \par 
We use exactly the same estimates as above for $\partial \eta
 _{j}$ and $\bar \partial \psi _{j}.$ $\hfill\blacksquare $\ \par 

\begin{Lmm}
We have\par 
\quad \quad \quad $\displaystyle A_{3}:=8\Re \lbrack \partial ((1-\left\vert{z}\right\vert
 ^{2})^{p+1})\bar \partial (\varphi _{A,j}(sz))\rbrack \leq $\par 
\quad \quad \quad \quad \quad \quad \quad $\displaystyle \lesssim \left\vert{\chi '}\right\vert (1-\left\vert{z}\right\vert
 ^{2})^{p}\psi _{j}^{q-1/2}+16q(p+1)(1-\left\vert{z}\right\vert
 ^{2})^{p}\psi _{j}^{q-1/2}\lesssim $\par 
\quad \quad \quad \quad \quad \quad \quad $\displaystyle \lesssim \left\vert{\chi '}\right\vert (1-\left\vert{z}\right\vert
 ^{2})^{p}d(sz,E)^{2q-1}+16q(p+1)(1-\left\vert{z}\right\vert
 ^{2})^{p}d(sz,E)^{2q-1}.$\par 
and\par 
\quad \quad \quad $\displaystyle -A_{3}\lesssim \left\vert{\chi '}\right\vert (1-\left\vert{z}\right\vert
 ^{2})^{p}\psi _{j}^{q-1/2}+8(p+1)(1-\left\vert{z}\right\vert
 ^{2})^{p-1/2}q\psi _{j}^{q}\lesssim $\par 
\quad \quad \quad \quad \quad $\displaystyle \lesssim \left\vert{\chi '}\right\vert (1-\left\vert{z}\right\vert
 ^{2})^{p}d(sz,E)^{2q-1}+8q(p+1)(1-\left\vert{z}\right\vert ^{2})^{p-1/2}d(sz,E)^{2q}.$
\end{Lmm}

      Proof.\ \par 
We have\ \par 
\quad \quad \quad $\displaystyle \bar \partial (\varphi _{A,j}(sz))=\psi _{j}^{q}\bar
 \partial \eta _{j}+\eta _{j}\bar \partial (\psi _{j}^{q})$\ \par 
For the term $\displaystyle \psi _{j}\bar \partial \eta _{j}$
 we proceed exactly as in lemma~\ref{CI26} to get\ \par 
\quad \quad \quad $\displaystyle \ \left\vert{\psi _{j}^{q}\bar \partial \eta _{j}}\right\vert
 \lesssim \left\vert{\chi '}\right\vert (1-\left\vert{z}\right\vert
 ^{2})^{p}\psi _{j}^{q-1/2}.$\ \par 
So it remains\ \par 
\quad \quad \quad $\displaystyle B:=8\eta _{j}\Re \lbrack \partial ((1-\left\vert{z}\right\vert
 ^{2})^{p+1})\bar \partial (\psi _{j}^{q})(sz))\rbrack =$\ \par 
\quad \quad \quad \quad \quad $\displaystyle =-8(p+1)(1-\left\vert{z}\right\vert ^{2})^{p}\eta
 _{j}(sz)\Re \lbrack \bar z\bar \partial (\psi _{j}^{q})(sz))\rbrack .$\ \par 
For this term we have by lemma~\ref{CI23}\ \par 
\quad \quad \quad $\displaystyle \bar \partial (\psi _{j})^{q}(z)=q\frac{(z-\alpha
 _{j})\left\vert{z-\alpha _{j}}\right\vert ^{2q-2}\left\vert{z-\beta
 _{j}}\right\vert ^{2q}}{\delta _{j}^{2q}}+q\frac{(z-\beta _{j})\left\vert{z-\alpha
 _{j}}\right\vert ^{2q}\left\vert{z-\beta _{j}}\right\vert ^{2q-2}}{\delta
 _{j}^{2q}}$\ \par 
hence $\displaystyle B=B_{1}+B_{2}$ with\ \par 
\quad \quad \quad $\displaystyle B_{1}:=-8(p+1)(1-\left\vert{z}\right\vert ^{2})^{p}\eta
 _{j}(sz)q\frac{\left\vert{sz-\alpha _{j}}\right\vert ^{2q-2}\left\vert{sz-\beta
 _{j}}\right\vert ^{2q}}{\delta _{j}^{2q}}\Re \lbrack \bar z(sz-\alpha
 _{j})\rbrack $\ \par 
and\ \par 
\quad \quad \quad $\displaystyle B_{2}:=-8(p+1)(1-\left\vert{z}\right\vert ^{2})^{p}\eta
 _{j}(sz)q\frac{\left\vert{sz-\alpha _{j}}\right\vert ^{2q}\left\vert{sz-\beta
 _{j}}\right\vert ^{2q-2}}{\delta _{j}^{2q}}\Re \lbrack \bar
 z(sz-\beta _{j})\rbrack .$\ \par 
\quad Now we shall apply lemma~\ref{CI6} to get that $\displaystyle
 \Re \lbrack \bar z(sz-\alpha _{j})\rbrack \leq 0$ iff $\displaystyle
 D(\frac{\alpha _{j}}{2},\ \frac{1}{2})$ so\ \par 
\quad \quad \quad $\displaystyle B_{1}\geq 0\iff z\in \Gamma _{j}\cap D(\frac{\alpha
 _{j}}{2},\ \frac{1}{2}).$\ \par 
The same way\ \par 
\quad \quad \quad $\displaystyle B_{2}\geq 0\iff z\in \Gamma _{j}\cap D(\frac{\beta
 _{j}}{2},\ \frac{1}{2}).$\ \par 
\quad If $\displaystyle z\notin D(\frac{\alpha _{j}}{2},\ \frac{1}{2}),$
 then we have that $\displaystyle (1-\left\vert{z}\right\vert
 ^{2})\leq 2\left\vert{z-\alpha _{j}}\right\vert ^{2}$ so we get\ \par 
\quad \quad \quad $\displaystyle \forall z\in \Gamma _{j}\cap D(\frac{\alpha _{j}}{2},\
 \frac{1}{2})^{c},\ -B_{1}\leq 8q(p+1)(1-\left\vert{z}\right\vert
 ^{2})^{p}\psi _{j}^{q}\left\vert{(sz-\alpha _{j})}\right\vert
 ^{-1}\leq $\ \par 
\quad \quad \quad \quad \quad \quad \quad \quad \quad \quad \quad \quad \quad \quad \quad \quad \quad $\displaystyle \leq 8q(p+1)(1-\left\vert{z}\right\vert ^{2})^{p-1/2}\psi
 _{j}^{q}.$\ \par 
The same way:\ \par 
\quad \quad \quad $\displaystyle \forall z\in \Gamma _{j}\cap D(\frac{\beta _{j}}{2},\
 \frac{1}{2})^{c},\ -B_{2}\leq 8q(p+1)(1-\left\vert{z}\right\vert
 ^{2})^{p}\psi _{j}^{q}\left\vert{(sz-\beta _{j})}\right\vert ^{-1}\leq $\ \par 
\quad \quad \quad \quad \quad \quad \quad \quad \quad \quad \quad \quad \quad \quad \quad \quad \quad $\displaystyle \leq 8q(p+1)(1-\left\vert{z}\right\vert ^{2})^{p-1/2}\psi
 _{j}^{q}.$\ \par 
Hence we get\ \par 
\quad \quad \quad $\displaystyle \forall z\in \Gamma _{j},\ -B\leq 16q(p+1)(1-\left\vert{z}\right\vert
 ^{2})^{p-1/2}\psi _{j}^{q}.$\ \par 
\quad Now we have $\displaystyle B_{1}\geq 0\iff z\in \Gamma _{j}\cap
 D(\frac{\alpha _{j}}{2},\ \frac{1}{2}),$ so\ \par 
\quad \quad \quad $\displaystyle B_{1}\leq 8(p+1)(1-\left\vert{z}\right\vert ^{2})^{p}\eta
 _{j}(sz)q\frac{\left\vert{sz-\alpha _{j}}\right\vert ^{2q-2}\left\vert{sz-\beta
 _{j}}\right\vert ^{2q}}{\delta _{j}^{2q}}\left\vert{sz-\alpha
 _{j}}\right\vert \leq $\ \par 
\quad \quad \quad \quad \quad \quad \quad $\displaystyle \leq 8(p+1)(1-\left\vert{z}\right\vert ^{2})^{p}q\left\vert{(sz-\alpha
 _{j})}\right\vert ^{-1}\psi _{j}^{q}.$\ \par 
And the same way\ \par 
\quad \quad \quad $\displaystyle B_{2}\leq 8(p+1)(1-\left\vert{z}\right\vert ^{2})^{p}\eta
 _{j}(sz)q\frac{\left\vert{sz-\beta _{j}}\right\vert ^{2q-2}\left\vert{sz-\alpha
 _{j}}\right\vert ^{2q}}{\delta _{j}^{2q}}\left\vert{sz-\beta
 _{j}}\right\vert \leq $\ \par 
\quad \quad \quad \quad \quad \quad \quad $\displaystyle \leq 8(p+1)(1-\left\vert{z}\right\vert ^{2})^{p}q\left\vert{(sz-\beta
 _{j})}\right\vert ^{-1}\psi _{j}^{q}.$\ \par 
Hence\ \par 
\quad \quad \quad $\displaystyle B\leq 8(p+1)(1-\left\vert{z}\right\vert ^{2})^{p}q(\left\vert{(sz-\alpha
 _{j})}\right\vert ^{-1}+\left\vert{(sz-\beta _{j})}\right\vert
 ^{-1})\psi _{j}^{q}.$\ \par 
\quad So we get\ \par 
\quad \quad \quad $\displaystyle \forall z\in \Gamma _{j},\ A_{3}:=8\Re \lbrack
 \partial ((1-\left\vert{z}\right\vert ^{2})^{p+1})\bar \partial
 (\varphi _{A,j}(sz))\rbrack \lesssim $\ \par 
\quad \quad \quad \quad \quad \quad \quad \quad \quad \quad \quad $\displaystyle \lesssim \left\vert{\chi '}\right\vert (1-\left\vert{z}\right\vert
 ^{2})^{p}\psi _{j}^{q-1/2}+16q(p+1)(1-\left\vert{z}\right\vert
 ^{2})^{p-1/2}\psi _{j}^{q}.$\ \par 
And\ \par 
\quad \quad \quad $\displaystyle -A_{3}\lesssim \left\vert{\chi '}\right\vert (1-\left\vert{z}\right\vert
 ^{2})^{p}\psi _{j}^{q-1/2}+8q(p+1)(1-\left\vert{z}\right\vert
 ^{2})^{p-1/2}\psi _{j}^{q}.$\ \par 
It remains to use lemma~\ref{CI25} to get the result. $\hfill\blacksquare
 $\ \par 
\ \par 
\quad We shall estimate $\displaystyle \triangle g_{A,s}(z)\log ^{+}\left\vert{f(sz)}\right\vert
 .$\ \par 

\begin{Prps}
~\label{eOP4}We have\par 
\quad \quad \quad $\displaystyle \triangle g_{A,s}(z)\lesssim 4p(p+1)(1-\left\vert{z}\right\vert
 ^{2})^{p-1}d(sz,E)^{2q}+(1-\left\vert{z}\right\vert ^{2})^{p}d(sz,E)^{2q-1}.$
\end{Prps}
\quad Proof.\ \par 
By use of $\displaystyle \triangle g_{A,s}(z)=A_{1}+A_{2}+A_{3}$
 and by the previous lemmas, we get for $\displaystyle z\in \Gamma _{j},$\ \par 
\quad \quad \quad $\displaystyle A_{1}=A'_{1}-A''_{1}\leq A'_{1}=4p(p+1)(1-\left\vert{z}\right\vert
 ^{2})^{p-1}\left\vert{z}\right\vert ^{2}\eta _{j}(sz)\psi _{j}^{q}(sz)\leq
 $\ \par 
\quad \quad \quad \quad \quad \quad \quad $\displaystyle \leq 4p(p+1)(1-\left\vert{z}\right\vert ^{2})^{p-1}d(sz,E)^{2q}.$\
 \par 
Then\ \par 
\quad \quad \quad $\displaystyle A_{2}=A_{2,1}+A_{2,2}+A_{2,3},$\ \par 
and, for $\displaystyle s\geq 1/2,$\ \par 
\quad \quad \quad $\displaystyle 0\leq A_{2,1}(s,z)\leq 4^{2q+2}q^{2}(1-\left\vert{z}\right\vert
 ^{2})^{p}c(\lambda )d(sz,E)^{2q-1}.$\ \par 
\quad \quad \quad $\displaystyle \ \left\vert{A_{2,2}(s,z)}\right\vert \lesssim
 (1-\left\vert{z}\right\vert ^{2})^{p}d(sz,E)^{2q-1}.$\ \par 
\quad \quad \quad $\displaystyle \ \left\vert{A_{2,3}}\right\vert \lesssim \left\vert{\chi
 '}\right\vert (1-\left\vert{z}\right\vert ^{2})^{p}d(sz,E)^{2q-1}.$\ \par 
Hence, for $\displaystyle \forall z\in \Gamma _{j},$\ \par 
\quad \quad \quad $\displaystyle A_{2}\lesssim (1-\left\vert{z}\right\vert ^{2})^{p}d(sz,E)^{2q-1}.$\
 \par 
Finally\ \par 
\quad \quad \quad $\displaystyle A_{3}\lesssim (1-\left\vert{z}\right\vert ^{2})^{p}d(sz,E)^{2q-1}+16q(p+1)(1-\left\vert{z}\right\vert
 ^{2})^{p-1/2}d(sz,E)^{2q}\lesssim (1-\left\vert{z}\right\vert
 ^{2})^{p}d(sz,E)^{2q-1}.$\ \par 
So we get\ \par 
\quad \quad \quad $\displaystyle \triangle g_{A,s}(z)\lesssim 4p(p+1)(1-\left\vert{z}\right\vert
 ^{2})^{p-1}d(sz,E)^{2q}+(1-\left\vert{z}\right\vert ^{2})^{p}d(sz,E)^{2q-1},$\
 \par 
which proves the proposition. $\hfill\blacksquare $\ \par 
\ \par 
\quad Now we shall estimate $\displaystyle -\triangle g_{A,s}(z)\log
 ^{-}\left\vert{f(sz)}\right\vert .$ We set:\ \par 
\quad \quad \quad $\displaystyle P_{{\mathbb{D}},A,-}(s):=\int_{{\mathbb{D}}}{(1-\left\vert{z}\right\vert
 ^{2})^{p-1}\left\vert{z}\right\vert ^{2}\varphi _{A}(sz)\log
 ^{-}\left\vert{fsz}\right\vert },$\ \par 
\quad \quad \quad $\displaystyle P_{{\mathbb{D}},A,+}(s):=\int_{{\mathbb{D}}}{(1-\left\vert{z}\right\vert
 ^{2})^{p-1}\varphi _{A}(sz)\log ^{+}\left\vert{fsz}\right\vert },$\ \par 
\quad \quad \quad $\displaystyle P_{-}(\delta ,u,s):=\int_{{\mathbb{D}}\backslash
 D(0,u)}{(1-\left\vert{z}\right\vert ^{2})^{p-1+\delta }\varphi
 _{A}(sz)\log ^{-}\left\vert{fsz}\right\vert }.$\ \par 
and\ \par 
\quad \quad \quad $\displaystyle P_{+}(\delta ,u,s):=\int_{D(0,u)}{(1-\left\vert{z}\right\vert
 ^{2})^{p-1+\delta }\varphi _{A}(sz)\log ^{+}\left\vert{fsz}\right\vert
 }.$\ \par 

\begin{Prps}
~\label{eOP5}We have, for $\displaystyle p>0,$\par 
\quad \quad \quad $\displaystyle -\int_{{\mathbb{D}}}{\triangle g_{A,s}(z)\log
 ^{-}\left\vert{f(sz)}\right\vert }\leq 2^{2q}P_{{\mathbb{D}},A,+}(s)+2{\times}4^{q}(1-u^{2})^{-2q}P_{+}(\frac{1}{2},u,s).$
\end{Prps}

      Proof.\ \par 
By use of $\displaystyle \triangle g_{A,s}(z)=A_{1}+A_{2}+A_{3}$
 and by the previous lemmas, we get for $\displaystyle \forall
 z\in \Gamma _{j},$\ \par 
\quad \quad \quad $\displaystyle -A_{1}=:-A'_{1}+A''_{1}=-4p(p+1)(1-\left\vert{z}\right\vert
 ^{2})^{p-1}\left\vert{z}\right\vert ^{2}\varphi _{A,j}(sz)+$\ \par 
\quad \quad \quad \quad \quad \quad \quad \quad \quad $\displaystyle +4(p+1)(1-\left\vert{z}\right\vert ^{2})^{p}d(sz,E)^{2q}.$\ \par 
Now\ \par 
\quad \quad \quad $\displaystyle -A_{2}=-A_{2,1}-A_{2,2}-A_{2,3}\leq -A_{2,2}-A_{2,3},$\ \par 
because $\displaystyle A_{2,1}\leq 0.$\ \par 
\quad We have\ \par 
\quad \quad \quad $\displaystyle \ \left\vert{A_{2,2}}\right\vert \lesssim (\left\vert{\chi
 '}\right\vert +\left\vert{\chi ''}\right\vert )(1-\left\vert{z}\right\vert
 ^{2})^{p}d(sz,E)^{2q-1}$\ \par 
and\ \par 
\quad \quad \quad $\displaystyle \ \left\vert{A_{2,3}}\right\vert \lesssim \left\vert{\chi
 '}\right\vert (1-\left\vert{z}\right\vert ^{2})^{p}d(sz,E)^{2q-1}$\ \par 
so\ \par 
\quad \quad \quad $\displaystyle -A_{2}\lesssim (\left\vert{\chi '}\right\vert
 +\left\vert{\chi ''}\right\vert )(1-\left\vert{z}\right\vert
 ^{2})^{p}d(sz,E)^{2q-1}.$\ \par 
Now\ \par 
\quad \quad \quad $\displaystyle -A_{3}\lesssim \left\vert{\chi '}\right\vert (1-\left\vert{z}\right\vert
 ^{2})^{p}d(sz,E)^{2q-1}+8q(p+1)(1-\left\vert{z}\right\vert ^{2})^{p-1/2}d(sz,E)^{2q}.$\
 \par 
\quad So grouping the terms we get\ \par 
\quad \quad \quad $\displaystyle -\triangle g_{A,s}(z)\lesssim (\left\vert{\chi
 '}\right\vert +\left\vert{\chi ''}\right\vert )(1-\left\vert{z}\right\vert
 ^{2})^{p}d(sz,E)^{2q-1}+8q(p+1)(1-\left\vert{z}\right\vert ^{2})^{p-1/2}d(sz,E)^{2q}-$\
 \par 
\quad \quad \quad \quad \quad \quad \quad \quad \quad $\displaystyle -4p(p+1)(1-\left\vert{z}\right\vert ^{2})^{p-1}\left\vert{z}\right\vert
 ^{2}\varphi _{A,j}(sz)+4(p+1)(1-\left\vert{z}\right\vert ^{2})^{p}\varphi
 _{A,j}(sz).$\ \par 
and\ \par 
\quad \quad \quad $\displaystyle -\triangle g_{A,s}(z)\log ^{-}\left\vert{f(sz)}\right\vert
 \lesssim (\left\vert{\chi '}\right\vert +\left\vert{\chi ''}\right\vert
 )(1-\left\vert{z}\right\vert ^{2})^{p}d(sz,E)^{2q-1}\log ^{-}\left\vert{f(sz)}\right\vert
 +$\ \par 
\quad \quad \quad \quad \quad \quad \quad \quad \quad \quad \quad \quad \quad \quad \quad \quad $\displaystyle +8q(p+1)(1-\left\vert{z}\right\vert ^{2})^{p-1/2}d(sz,E)^{2q}\log
 ^{-}\left\vert{f(sz)}\right\vert -$\ \par 
\quad \quad \quad \quad \quad \quad \quad \quad \quad \quad \quad \quad \quad \quad \quad \quad $\displaystyle -4p(p+1)(1-\left\vert{z}\right\vert ^{2})^{p-1}\left\vert{z}\right\vert
 ^{2}\varphi _{A,j}(sz)\log ^{-}\left\vert{f(sz)}\right\vert +$\ \par 
\quad \quad \quad \quad \quad \quad \quad \quad \quad \quad \quad \quad \quad \quad \quad \quad $\displaystyle +4(p+1)(1-\left\vert{z}\right\vert ^{2})^{p}\varphi
 _{A,j}(sz)\log ^{-}\left\vert{f(sz)}\right\vert .$\ \par 
\quad For the first term, because on $\displaystyle (\left\vert{\chi
 '}\right\vert +\left\vert{\chi ''}\right\vert )\neq 0$ we have
 $\displaystyle d(sz,E)\leq 3(1-\left\vert{sz}\right\vert ^{2}),$ we get\ \par 
\quad \quad \quad $\displaystyle B_{1}:=(\left\vert{\chi '}\right\vert +\left\vert{\chi
 ''}\right\vert )(1-\left\vert{z}\right\vert ^{2})^{p}d(sz,E)^{2q-1}\log
 ^{-}\left\vert{f(sz)}\right\vert \lesssim $\ \par 
\quad \quad \quad \quad \quad \quad \quad \quad \quad \quad \quad \quad \quad $\displaystyle \lesssim (1-\left\vert{z}\right\vert ^{2})^{p}(1-\left\vert{sz}\right\vert
 ^{2})^{2q-1}\log ^{-}\left\vert{f(sz)}\right\vert $\ \par 
hence, passing in polar coordinates,\ \par 
\quad \quad \quad $\displaystyle \ \int_{{\mathbb{D}}}{B_{1}}\lesssim \int_{0}^{1}{(1-\rho
 ^{2})^{p}(1-s^{2}\rho ^{2})^{2q-1}\lbrace \int_{{\mathbb{T}}}{\log
 ^{-}\left\vert{f(s\rho e^{i\theta })}\right\vert }\rbrace \rho d\rho }.$\ \par 
By the subharmonicity of $\displaystyle \log \left\vert{f(sz)}\right\vert
 $ and $\displaystyle \ \left\vert{f(0)}\right\vert =1,$ we get\ \par 
\quad \quad \quad $\displaystyle \ \int_{{\mathbb{T}}}{\log ^{-}\left\vert{f(s\rho
 e^{i\theta })}\right\vert }\leq \int_{{\mathbb{T}}}{\log ^{+}\left\vert{f(s\rho
 e^{i\theta })}\right\vert },$\ \par 
hence\ \par 
\quad \quad \quad $\displaystyle \ \int_{{\mathbb{D}}}{B_{1}}\lesssim \int_{{\mathbb{D}}}{(1-\left\vert{z}\right\vert
 ^{2})^{p}(1-\left\vert{sz}\right\vert ^{2})^{2q-1}\log ^{+}\left\vert{f(sz)}\right\vert
 }.$\ \par 
\quad For $\displaystyle B_{2}:=8q(p+1)(1-\left\vert{z}\right\vert
 ^{2})^{p-1/2}d(sz,E)^{2q}\log ^{-}\left\vert{f(sz)}\right\vert
 ,$ we use the substitution lemma~\ref{3_CIZ5} with $\delta =1/2,$
 to get:\ \par 
\quad \quad \quad $\displaystyle \ \int_{{\mathbb{D}}}{B_{2}}\leq 4^{q}(1-u^{2})^{-2q}P_{+}(\frac{1}{2},u)+(1-u^{2})^{1/4}u^{-2}P_{-}(\frac{1}{4},u,s).$\
 \par 
\quad For $\displaystyle p>0,$ because $\displaystyle \varphi _{A}(z)\lesssim
 d(z,E)^{2q},$ we get:\ \par 
\quad \quad \quad $\displaystyle \ \int_{{\mathbb{D}}}{B_{2}}\lesssim 4^{q}(1-u^{2})^{-2q}\int_{D(0,u)}{(1-\left\vert{z}\right\vert
 ^{2})^{p-1/2}d(sz,E)^{2q}\log ^{+}\left\vert{fsz}\right\vert }+$\ \par 
\quad \quad \quad \quad \quad \quad \quad \quad \quad $\displaystyle +(1-u^{2})^{1/4}u^{-2}\int_{{\mathbb{D}}}{(1-\left\vert{z}\right\vert
 ^{2})^{p-3/4}\left\vert{z}\right\vert ^{2}\varphi _{A}(sz)\log
 ^{-}\left\vert{fsz}\right\vert }.$\ \par 
The same for the last term with $\delta =1$ and we get that\ \par 
\quad \quad \quad $\displaystyle B_{4}:=4(p+1)(1-\left\vert{z}\right\vert ^{2})^{p}d(sz,E)^{2q}\log
 ^{-}\left\vert{f(sz)}\right\vert $\ \par 
verifies\ \par 
\quad \quad \quad $\displaystyle \ \int_{{\mathbb{D}}}{B_{4}}\lesssim 4^{q}(1-u^{2})^{-2q}P_{+}(1,u,s)+(1-u^{2})^{1/2}u^{-2}P_{-}(\frac{1}{2},u,s).$\
 \par 
\quad Now it remains the "good" term\ \par 
\quad \quad \quad $\displaystyle B_{3}:=-4p(p+1)(1-\left\vert{z}\right\vert ^{2})^{p-1}\left\vert{z}\right\vert
 ^{2}\varphi _{A,j}(sz)\log ^{-}\left\vert{f(sz)}\right\vert $\ \par 
and, if $\displaystyle p>0,$ we choose $\displaystyle 1-u^{2}$
 small enough to get that\ \par 
\quad \quad \quad $\displaystyle (1-u^{2})^{1/4}u^{-2}\int_{{\mathbb{D}}}{(1-\left\vert{z}\right\vert
 ^{2})^{p-3/4}\left\vert{z}\right\vert ^{2}\varphi _{A}(sz)\log
 ^{-}\left\vert{fsz}\right\vert }+$\ \par 
\quad \quad \quad \quad \quad $\displaystyle +(1-u^{2})^{1/2}u^{-2}\int_{{\mathbb{D}}}{(1-\left\vert{z}\right\vert
 ^{2})^{p-1/2}\left\vert{z}\right\vert ^{2}\varphi _{A}(sz)\log
 ^{-}\left\vert{fsz}\right\vert }-$\ \par 
\quad \quad \quad \quad \quad \quad \quad $\displaystyle -4p(p+1)\int_{{\mathbb{D}}}{(1-\left\vert{z}\right\vert
 ^{2})^{p-1}\left\vert{z}\right\vert ^{2}\varphi _{A}(sz)\log
 ^{-}\left\vert{f(sz)}\right\vert }\leq 0.$\ \par 
So it remains:\ \par 
\quad \quad \quad $\displaystyle -\int_{{\mathbb{D}}}{\triangle g_{A,s}(z)\log
 ^{-}\left\vert{f(sz)}\right\vert }\leq \int_{{\mathbb{D}}}{(1-\left\vert{z}\right\vert
 ^{2})^{p}(1-\left\vert{sz}\right\vert ^{2})^{2q-1}\log ^{+}\left\vert{f(sz)}\right\vert
 }+$\ \par 
\quad \quad \quad \quad \quad \quad \quad \quad \quad \quad \quad \quad \quad \quad \quad \quad \quad \quad \quad $\displaystyle +2^{q}(1-u^{2})^{-2q}P_{+}(\frac{1}{2},u,s)+$\ \par 
\quad \quad \quad \quad \quad \quad \quad \quad \quad \quad \quad \quad \quad \quad \quad \quad \quad \quad \quad \quad $\displaystyle +2^{q}(1-u^{2})^{-2q}P_{+}(1,u,s)\leq $\ \par 
\quad \quad \quad \quad \quad \quad \quad \quad \quad \quad \quad \quad \quad \quad \quad $\displaystyle \leq \int_{{\mathbb{D}}}{(1-\left\vert{z}\right\vert
 ^{2})^{p}(1-\left\vert{sz}\right\vert ^{2})^{2q-1}\log ^{+}\left\vert{f(sz)}\right\vert
 }+$\ \par 
\quad \quad \quad \quad \quad \quad \quad \quad \quad \quad \quad \quad \quad \quad \quad \quad \quad \quad \quad $\displaystyle +2{\times}4^{q}(1-u^{2})^{-2q}P_{+}(\frac{1}{2},u,s).$\ \par 
Now $\displaystyle (1-\left\vert{z}\right\vert ^{2})\leq (1-\left\vert{sz}\right\vert
 ^{2}),$ and $\displaystyle (1-\left\vert{sz}\right\vert ^{2})^{2q}\leq
 2^{2q}\varphi _{A}(sz)$ so we get\ \par 
\quad \quad \quad $\displaystyle \ \int_{{\mathbb{D}}}{(1-\left\vert{z}\right\vert
 ^{2})^{p}(1-\left\vert{sz}\right\vert ^{2})^{2q-1}\log ^{+}\left\vert{f(sz)}\right\vert
 }\leq 2^{2q}\int_{{\mathbb{D}}}{(1-\left\vert{z}\right\vert
 ^{2})^{p-1}\varphi _{A}(sz)\log ^{+}\left\vert{f(sz)}\right\vert },$\ \par 
so putting it, we get\ \par 
\quad \quad \quad $\displaystyle -\int_{{\mathbb{D}}}{\triangle g_{A,s}(z)\log
 ^{-}\left\vert{f(sz)}\right\vert }\leq 2^{2q}\int_{{\mathbb{D}}}{(1-\left\vert{z}\right\vert
 ^{2})^{p-1}\varphi _{A}(sz)\log ^{+}\left\vert{f(sz)}\right\vert }+$\ \par 
\quad \quad \quad \quad \quad \quad \quad \quad \quad \quad \quad \quad \quad \quad \quad $\displaystyle +2{\times}4^{q}(1-u^{2})^{-2q}P_{+}(\frac{1}{2},u),$\ \par 
which ends the proof. $\hfill\blacksquare $\ \par 

\begin{Prps}
~\label{eOP7}We have, for $\displaystyle p=0,$\par 
\quad \quad \quad $\displaystyle -\int_{{\mathbb{D}}}{\triangle g_{A,s}(z)\log
 ^{-}\left\vert{f(sz)}\right\vert \lesssim }\int_{{\mathbb{D}}}{(1-\left\vert{sz}\right\vert
 ^{2})^{2q-1}\log ^{+}\left\vert{f(sz)}\right\vert }+$\par 
\quad \quad \quad \quad \quad \quad \quad \quad \quad $\displaystyle +2{\times}4^{q}(1-u^{2})^{-2q}\int_{D(0,u)}{(1-\left\vert{z}\right\vert
 ^{2})^{-1/2}d(sz,E)^{2q}\log ^{+}\left\vert{fsz}\right\vert }+$\par 
\quad \quad \quad \quad \quad \quad \quad \quad \quad $\displaystyle +2(1-u^{2})^{1/4}u^{-2}\int_{{\mathbb{D}}\backslash
 D(0,u)}{(1-\left\vert{z}\right\vert ^{2})^{-3/4}\left\vert{z}\right\vert
 ^{2}\varphi _{A}(sz)\log ^{-}\left\vert{fsz}\right\vert }.$
\end{Prps}
\quad Proof.\ \par 
\quad For $\displaystyle p=0,$ there is no "good" term and we have\ \par 
\quad \quad \quad $\displaystyle \ \int_{{\mathbb{D}}}{B_{1}}\lesssim \int_{{\mathbb{D}}}{(1-\left\vert{sz}\right\vert
 ^{2})^{2q-1}\log ^{+}\left\vert{f(sz)}\right\vert }.$\ \par 
And\ \par 
\quad \quad \quad $\displaystyle \ \int_{{\mathbb{D}}}{B_{2}}\lesssim 4^{q}(1-u^{2})^{-2q}\int_{D(0,u)}{(1-\left\vert{z}\right\vert
 ^{2})^{-1/2}d(sz,E)^{2q}\log ^{+}\left\vert{fsz}\right\vert }+$\ \par 
\quad \quad \quad \quad \quad \quad \quad \quad \quad $\displaystyle +(1-u^{2})^{1/4}u^{-2}\int_{{\mathbb{D}}\backslash
 D(0,u)}{(1-\left\vert{z}\right\vert ^{2})^{-3/4}\left\vert{z}\right\vert
 ^{2}\varphi _{A}(sz)\log ^{-}\left\vert{fsz}\right\vert }.$\ \par 
The same for the last term with $\delta =1$ and we get\ \par 
\quad \quad \quad $\displaystyle B_{4}:=4d(sz,E)^{2q}\log ^{-}\left\vert{f(sz)}\right\vert
 $\ \par 
verifies\ \par 
\quad \quad \quad $\displaystyle \ \int_{{\mathbb{D}}}{B_{4}}\lesssim 4^{q}(1-u^{2})^{-2q}P_{+}(1,u)+(1-u^{2})^{1/2}u^{-2}P_{-}(\frac{1}{2},u,s).$\
 \par 
So adding:\ \par 
\quad \quad \quad $\displaystyle -\int_{{\mathbb{D}}}{\triangle g_{A,s}(z)\log
 ^{-}\left\vert{f(sz)}\right\vert \lesssim }\int_{{\mathbb{D}}}{(1-\left\vert{sz}\right\vert
 ^{2})^{2q-1}\log ^{+}\left\vert{f(sz)}\right\vert }+$\ \par 
\quad \quad \quad \quad \quad \quad \quad \quad \quad $\displaystyle +4^{q}(1-u^{2})^{-2q}\int_{D(0,u)}{(1-\left\vert{z}\right\vert
 ^{2})^{-1/2}d(sz,E)^{2q}\log ^{+}\left\vert{fsz}\right\vert }+$\ \par 
\quad \quad \quad \quad \quad \quad \quad \quad \quad $\displaystyle +4^{q}(1-u^{2})^{-2q}\int_{D(0,u)}{d(sz,E)^{2q}\log
 ^{+}\left\vert{fsz}\right\vert }+$\ \par 
\quad \quad \quad \quad \quad \quad \quad \quad \quad $\displaystyle +(1-u^{2})^{1/4}u^{-2}\int_{{\mathbb{D}}\backslash
 D(0,u)}{(1-\left\vert{z}\right\vert ^{2})^{-3/4}\left\vert{z}\right\vert
 ^{2}\varphi _{A}(sz)\log ^{-}\left\vert{fsz}\right\vert }+$\ \par 
\quad \quad \quad \quad \quad \quad \quad \quad \quad $\displaystyle +(1-u^{2})^{1/2}u^{-2}\int_{{\mathbb{D}}\backslash
 D(0,u)}{(1-\left\vert{z}\right\vert ^{2})^{-1/2}\left\vert{z}\right\vert
 ^{2}\varphi _{A}(sz)\log ^{-}\left\vert{fsz}\right\vert }.$\ \par 
And\ \par 
\quad \quad \quad $\displaystyle -\int_{{\mathbb{D}}}{\triangle g_{A,s}(z)\log
 ^{-}\left\vert{f(sz)}\right\vert \lesssim }\int_{{\mathbb{D}}}{(1-\left\vert{sz}\right\vert
 ^{2})^{2q-1}\log ^{+}\left\vert{f(sz)}\right\vert }+$\ \par 
\quad \quad \quad \quad \quad \quad \quad \quad \quad $\displaystyle +2{\times}4^{q}(1-u^{2})^{-2q}\int_{D(0,u)}{(1-\left\vert{z}\right\vert
 ^{2})^{-1/2}d(sz,E)^{2q}\log ^{+}\left\vert{fsz}\right\vert }+$\ \par 
\quad \quad \quad \quad \quad \quad \quad \quad \quad $\displaystyle +2(1-u^{2})^{1/4}u^{-2}\int_{{\mathbb{D}}\backslash
 D(0,u)}{(1-\left\vert{z}\right\vert ^{2})^{-3/4}\left\vert{z}\right\vert
 ^{2}\varphi _{A}(sz)\log ^{-}\left\vert{fsz}\right\vert },$\ \par 
which ends the proof. $\hfill\blacksquare $\ \par 

\begin{Prps}
~\label{eOP6}We have, for $\displaystyle p>0,$\par 
\quad \quad \quad $\displaystyle \ \int_{{\mathbb{D}}}{\triangle g_{A,s}(z)\log
 \left\vert{f(sz)}\right\vert }\lesssim \lbrack 2^{2q}+4p(p+1)+2\rbrack
 P_{{\mathbb{D}},A,+}(s)+2{\times}4^{q}(1-u^{2})^{-2q}P_{+}(\frac{1}{2},u,s),$\par
 
\quad Proof.
\end{Prps}
From\ \par 
\quad \quad \quad $\displaystyle \triangle g_{A,s}(z)\log \left\vert{f(sz)}\right\vert
 =\triangle g_{A,s}(z)\log ^{+}\left\vert{f(sz)}\right\vert -\triangle
 g_{A,s}(z)\log ^{-}\left\vert{f(sz)}\right\vert $\ \par 
using proposition~\ref{eOP4} we have, using $\displaystyle (1-\left\vert{z}\right\vert
 ^{2})\leq 2d(sz,E),$\ \par 
\quad \quad \quad $\displaystyle \triangle g_{A,s}(z)\log ^{+}\left\vert{f(sz)}\right\vert
 \lesssim 4p(p+1)(1-\left\vert{z}\right\vert ^{2})^{p-1}d(sz,E)^{2q}\log
 ^{+}\left\vert{fsz}\right\vert +$\ \par 
\quad \quad \quad \quad \quad \quad \quad \quad \quad \quad \quad \quad \quad $\displaystyle +(1-\left\vert{z}\right\vert ^{2})^{p}d(sz,E)^{2q-1}\log
 ^{+}\left\vert{fsz}\right\vert \leq $\ \par 
\quad \quad \quad \quad \quad \quad \quad \quad \quad \quad \quad $\displaystyle \leq \lbrack 4p(p+1)+2\rbrack (1-\left\vert{z}\right\vert
 ^{2})^{p-1}d(sz,E)^{2q}\log ^{+}\left\vert{fsz}\right\vert .$\ \par 
And using proposition~\ref{eOP5} we have\ \par 
\quad \quad \quad $\displaystyle -\int_{{\mathbb{D}}}{\triangle g_{A,s}(z)\log
 ^{-}\left\vert{f(sz)}\right\vert }\leq 2^{2q}P_{{\mathbb{D}},A,+}(s)+2{\times}4^{q}(1-u^{2})^{-2q}P_{+}(\frac{1}{2},u,s).$\
 \par 
Hence\ \par 
\quad \quad \quad $\displaystyle \ \int_{{\mathbb{D}}}{\triangle g_{A,s}(z)\log
 \left\vert{f(sz)}\right\vert }\lesssim 2^{2q}P_{{\mathbb{D}},A,+}(s)+2{\times}4^{q}(1-u^{2})^{-2q}P_{+}(\frac{1}{2},u,s)+$\
 \par 
\quad \quad \quad \quad \quad \quad \quad \quad \quad \quad \quad \quad \quad \quad \quad $\displaystyle +\lbrack 4p(p+1)+2\rbrack \int_{{\mathbb{D}}}{(1-\left\vert{z}\right\vert
 ^{2})^{p-1}d(sz,E)^{2q}\log ^{+}\left\vert{fsz}\right\vert },$\ \par 
so\ \par 
\quad \quad \quad $\displaystyle \ \int_{{\mathbb{D}}}{\triangle g_{A,s}(z)\log
 \left\vert{f(sz)}\right\vert }\lesssim \lbrack 2^{2q}+4p(p+1)+2\rbrack
 P_{{\mathbb{D}},A,+}(s)+2{\times}4^{q}(1-u^{2})^{-2q}P_{+}(\frac{1}{2},u,s),$\
 \par 
which ends the proof. $\hfill\blacksquare $\ \par 

\section{The case $p>0.$}
\quad Recall that\ \par 
\quad \quad \quad $\displaystyle \forall z\in \Gamma _{j},\ \varphi _{j}(z):=\eta
 _{j}(z)\psi _{j}(z)^{q}+(1-\left\vert{z}\right\vert ^{2})^{2q},\
 \forall z\in \Gamma _{E},\ \varphi _{E}(z):=(1-\left\vert{z}\right\vert
 ^{2})^{2q},$\ \par 
and by lemma~\ref{eOPC2} we have that there is a function $\displaystyle
 \varphi \in {\mathcal{C}}^{\infty }({\mathbb{D}})$ such that
 $\varphi $ coincides with $\displaystyle \varphi _{j}$ and $\displaystyle
 \varphi _{E}$ in their domains of definition. Moreover we have
 for $\displaystyle 0\leq s<1$ and $\displaystyle q>0,\ g_{s}(z):=(1-\left\vert{z}\right\vert
 ^{2})^{p+1}\varphi (sz)\in {\mathcal{C}}^{\infty }(\bar {\mathbb{D}})$
 so we can apply the Green formula to it. Recall that $\displaystyle
 f_{s}(z):=f(sz).$\ \par 
\quad With the "zero" formula: $\Delta \log \left\vert{f_{s}}\right\vert
 =\sum_{a\in Z(f_{s})}{\delta _{a}}$ we get\ \par 
\quad \quad \quad $\displaystyle \ \sum_{a\in Z(f_{s})}{g_{s}(a)}=\int_{{\mathbb{D}}}{\log
 \left\vert{f(sz)}\right\vert \triangle g_{s}(z)}+\int_{{\mathbb{T}}}{(g_{s}\partial
 _{n}\log \left\vert{f(sz)}\right\vert -\log \left\vert{f(sz)}\right\vert
 \partial _{n}g_{s})}.$\ \par 
So, because $\displaystyle g_{s}=0$ on $\displaystyle {\mathbb{T}},$\ \par 
\quad \quad \quad $\displaystyle \ \sum_{a\in Z(f_{s})}{g_{s}(a)}=\int_{{\mathbb{D}}}{\log
 \left\vert{f(sz)}\right\vert \triangle g_{s}(z)}-\int_{{\mathbb{T}}}{\log
 \left\vert{f(se^{i\theta })}\right\vert \partial _{n}g_{s}(e^{i\theta
 })}.$\ \par 
If, moreover $\displaystyle p>0,\ \partial _{n}g_{s}=0$ on $\displaystyle
 {\mathbb{T}},$ hence $\displaystyle \ \sum_{a\in Z(f_{s})}{g_{s}(a)}=\int_{{\mathbb{D}}}{\log
 \left\vert{f(sz)}\right\vert \triangle g_{s}(z)}.$\ \par 
\ \par 
\quad So we have to compute\ \par 
\quad \quad \quad $\displaystyle \triangle g_{s}(z)\log \left\vert{f(sz)}\right\vert
 =\triangle g_{s}(z)\log ^{+}\left\vert{f(sz)}\right\vert -\triangle
 g_{s}(z)\log ^{-}\left\vert{f(sz)}\right\vert .$\ \par 
We have $\Delta g_{s}=4\bar \partial \partial g_{s}$ hence\ \par 
\quad \quad \quad $\displaystyle \Delta g_{s}(z)=\Delta \lbrack (1-\left\vert{z}\right\vert
 ^{2})^{p+1}\varphi (sz)\rbrack =\varphi (sz)\Delta \lbrack (1-\left\vert{z}\right\vert
 ^{2})^{p+1}\rbrack +(1-\left\vert{z}\right\vert ^{2})^{p+1}\Delta
 \lbrack \varphi (sz)\rbrack +$\ \par 
\quad \quad \quad \quad \quad \quad \quad $\displaystyle +8\Re \lbrack \partial ((1-\left\vert{z}\right\vert
 ^{2})^{p+1})\bar \partial (\varphi (sz))\rbrack .$\ \par 
\ \par 
\quad Recall that, with $\displaystyle \varphi _{A,j}(z):=\eta _{j}(z)\psi
 _{j}(z)^{q}$ and $\displaystyle \varphi _{C,j}(z):=(1-\left\vert{z}\right\vert
 ^{2})^{2q},$ we have\ \par 
\quad \quad \quad $\displaystyle \forall z\in \Gamma _{j},\ \varphi _{j}(z):=\varphi
 _{A,j}(z)+\varphi _{C,j}(z),$\ \par 
and\ \par 
\quad \quad \quad $\displaystyle g_{C,s}(z):=(1-\left\vert{z}\right\vert ^{2})^{p+1}\varphi
 _{C}(sz)=(1-\left\vert{z}\right\vert ^{2})^{p+1}(1-\left\vert{sz}\right\vert
 ^{2})^{2q},$\ \par 
and\ \par 
\quad \quad \quad $\displaystyle g_{A,s}(z):=(1-\left\vert{z}\right\vert ^{2})^{p+1}\sum_{j\in
 {\mathbb{N}}}{{\11}_{\Gamma _{j}}(z)\varphi _{A,j}(sz)}.$\ \par 
Now we are in position to apply the previous results. By proposition~\ref{CI7}
 we get:\ \par 
\quad \quad \quad $\displaystyle \ \int_{{\mathbb{D}}}{\Delta g_{C,s}(z)\log ^{+}\left\vert{f(sz)}\right\vert
 }\leq c(p,q)\int_{{\mathbb{D}}}{(1-\left\vert{z}\right\vert
 ^{2})^{p-1}\varphi _{C}(sz)\log ^{+}\left\vert{f(sz)}\right\vert }.$\ \par 
And by proposition~\ref{CI4} we get:\ \par 
\quad \quad \quad $\displaystyle -\int_{{\mathbb{D}}}{\Delta g_{C,s}(z)\log ^{-}\left\vert{f(sz)}\right\vert
 }\leq 4\lbrack (p+1)+2sq\rbrack \int_{{\mathbb{D}}}{(1-\left\vert{z}\right\vert
 ^{2})^{p}(1-\left\vert{sz}\right\vert ^{2})^{2q}\log ^{+}\left\vert{f(sz)}\right\vert
 }.$\ \par 
\quad So adding:\ \par 
\quad \quad \quad $\displaystyle \ \int_{{\mathbb{D}}}{\Delta g_{C,s}(z)\log \left\vert{f(sz)}\right\vert
 }\leq c(p,q)\int_{{\mathbb{D}}}{(1-\left\vert{z}\right\vert
 ^{2})^{p-1}\varphi _{C}(sz)\log ^{+}\left\vert{f(sz)}\right\vert }+$\ \par 
\quad \quad \quad \quad \quad \quad \quad \quad \quad \quad \quad \quad \quad \quad \quad $\displaystyle +4\lbrack (p+1)+2sq\rbrack \int_{{\mathbb{D}}}{(1-\left\vert{z}\right\vert
 ^{2})^{p}(1-\left\vert{sz}\right\vert ^{2})^{2q}\log ^{+}\left\vert{f(sz)}\right\vert
 }.$\ \par 
\quad Now by proposition~\ref{eOP6} we get:\ \par 
\quad \quad \quad $\displaystyle \ \int_{{\mathbb{D}}}{\triangle g_{A,s}(z)\log
 \left\vert{f(sz)}\right\vert }\lesssim \lbrack 2^{2q}+4p(p+1)+2\rbrack
 \int_{{\mathbb{D}}}{(1-\left\vert{z}\right\vert ^{2})^{p-1}\varphi
 _{A}(sz)\log ^{+}\left\vert{fsz}\right\vert }+$\ \par 
\quad \quad \quad \quad \quad \quad \quad \quad \quad $\displaystyle +2{\times}4^{q}(1-u^{2})^{-2q}\int_{D(0,u)}{(1-\left\vert{z}\right\vert
 ^{2})^{p-1/2}\varphi _{A}(sz)\log ^{+}\left\vert{fsz}\right\vert }.$\ \par 
Adding, because $\varphi =\varphi _{A}+\varphi _{C},$ we get\ \par 

\begin{Thrm}
~\label{eP0} We have:\par 
\quad \quad \quad $\displaystyle \ \int_{{\mathbb{D}}}{\triangle g_{s}(z)\log \left\vert{f(sz)}\right\vert
 }\lesssim \int_{{\mathbb{D}}}{(1-\left\vert{z}\right\vert ^{2})^{p-1}\varphi
 (sz)\log ^{+}\left\vert{fsz}\right\vert }.$
\end{Thrm}
\quad Proof.\ \par 
This is clear because, in the second term:\ \par 
\quad \quad \quad $\displaystyle (1-\left\vert{z}\right\vert ^{2})^{p-1/2}\varphi
 (sz)\log ^{+}\left\vert{fsz}\right\vert \leq (1-\left\vert{z}\right\vert
 ^{2})^{p-1}\varphi (sz)\log ^{+}\left\vert{fsz}\right\vert .$
 $\hfill\blacksquare $\ \par 
\ \par 
So we are lead to\ \par 

\begin{Dfnt}
Let $\displaystyle E=\bar E\subset {\mathbb{T}}.$ We say that
 an holomorphic function $f$ is in the generalised Nevanlinna
 class $\displaystyle {\mathcal{N}}_{\varphi ,p}({\mathbb{D}})$
 for $\displaystyle p>0$ if  $\displaystyle \exists \delta >0,\
 \delta <1$ such that\par 
\quad \quad \quad $\displaystyle \ {\left\Vert{f}\right\Vert}_{{\mathcal{N}}_{\varphi
 ,p}}:=\sup _{1-\delta <s<1}\int_{{\mathbb{D}}}{(1-\left\vert{z}\right\vert
 )^{p-1}\varphi (sz)\log ^{+}\left\vert{f(sz)}\right\vert }<\infty .$
\end{Dfnt}
\quad And we proved the Blaschke type condition:\ \par 

\begin{Thrm}
~\label{eP1}Let $\displaystyle E=\bar E\subset {\mathbb{T}}.$
 Suppose $\displaystyle q>0$ and $\displaystyle f\in {\mathcal{N}}_{\varphi
 ,p}({\mathbb{D}})$ with $\displaystyle \ \left\vert{f(0)}\right\vert
 =1,$ then\par 
\quad \quad \quad $\displaystyle \ \sum_{a\in Z(f)}{(1-\left\vert{a}\right\vert
 ^{2})^{1+p}\varphi (a)}\leq c(\varphi ){\left\Vert{f}\right\Vert}_{{\mathcal{N}}_{\varphi
 ,p}}.$
\end{Thrm}

\begin{Crll}
~\label{eP2}Let $\displaystyle E=\bar E\subset {\mathbb{T}}.$
 Suppose $\displaystyle q\in {\mathbb{R}}$ and $\displaystyle
 f\in {\mathcal{N}}_{d(\cdot ,E)^{q},p}({\mathbb{D}})$ with $\displaystyle
 \ \left\vert{f(0)}\right\vert =1,$ then\par 
\quad \quad \quad $\displaystyle \ \sum_{a\in Z(f)}{(1-\left\vert{a}\right\vert
 ^{2})^{1+p}d(a,E)^{q}}\leq c(\varphi ){\left\Vert{f}\right\Vert}_{{\mathcal{N}}_{d(\cdot
 ,E)^{q},p}}.$
\end{Crll}
\quad Proof.\ \par 
By use of lemma~\ref{eOPC3}, we have\ \par 
\quad \quad \quad $\displaystyle \ \sum_{a\in Z(f_{s})}{g_{s}(a)}=\int_{{\mathbb{D}}}{\log
 \left\vert{f(sz)}\right\vert \triangle g_{s}(z)},$\ \par 
hence, by theorem~\ref{eP0},\ \par 
\quad \quad \quad $\displaystyle \ \sum_{a\in Z(f_{s})}{g_{s}(a)}\lesssim \int_{{\mathbb{D}}}{(1-\left\vert{z}\right\vert
 ^{2})^{p-1}\varphi (sz)\log ^{+}\left\vert{fsz}\right\vert },$\ \par 
the constant in $\displaystyle \lesssim $ being independent of
 $\displaystyle s<1.$ It remains to apply lemma~\ref{6_A0} to
 get that, for any $\displaystyle 1>\delta >0$ we have:\ \par 
\quad \quad \quad $\displaystyle \ \sum_{a\in Z(f)}{(1-\left\vert{a}\right\vert
 ^{2})^{p+1}\varphi (a)}\leq \sup _{1-\delta <s<1}\int_{{\mathbb{D}}}{(1-\left\vert{z}\right\vert
 ^{2})^{p-1}\varphi (sz)\log ^{+}\left\vert{f(sz)}\right\vert },$\ \par 
which ends the proof of theorem~\ref{eP1}. $\hfill\blacksquare $\ \par 
\quad To prove corollary~\ref{eP2}, we use lemma~\ref{CI25} and lemma~\ref{eOPC2}
 which give that $\varphi (z)$ is equivalent to $\displaystyle
 d(z,E)^{2q}.$ $\hfill\blacksquare $\ \par 

\section{The case $p=0.$}
\quad This time we have $\displaystyle g_{s}(z):=(1-\left\vert{z}\right\vert
 ^{2})\varphi (sz)\in {\mathcal{C}}^{\infty }(\bar {\mathbb{D}})$ hence\ \par 
\quad \quad \quad $\displaystyle \ \sum_{a\in Z(f_{s})}{g_{s}(a)}=\int_{{\mathbb{D}}}{\log
 \left\vert{f(sz)}\right\vert \triangle g_{s}(z)}-\int_{{\mathbb{T}}}{\log
 \left\vert{f(se^{i\theta })}\right\vert \partial _{n}g_{s}(e^{i\theta
 })},$\ \par 
with\ \par 
\quad \quad \quad $\displaystyle \partial _{n}g_{s}(z)=-2\varphi (sz)+(1-\left\vert{z}\right\vert
 ^{2})\partial _{n}\varphi (sz),$\ \par 
so\ \par 
\quad \quad \quad $\displaystyle \forall e^{i\theta }\in {\mathbb{T}},\ \partial
 _{n}g_{s}(e^{i\theta })=-2\varphi (se^{i\theta })$\ \par 
hence\ \par 
\quad \quad \quad \begin{equation}  \ \sum_{a\in Z(f_{s})}{g_{s}(a)}=\int_{{\mathbb{D}}}{\log
 \left\vert{f(sz)}\right\vert \triangle g_{s}(z)}+2\int_{{\mathbb{T}}}{\varphi
 (e^{i\theta })\log \left\vert{f(se^{i\theta })}\right\vert }.\label{eP3}\end{equation}\
 \par 
\quad Now we are in position to apply the previous results. By proposition~\ref{CI7},
 we get:\ \par 
\quad \quad \quad $\displaystyle \ \int_{{\mathbb{D}}}{\Delta g_{C,s}(z)\log ^{+}\left\vert{f(sz)}\right\vert
 }\leq c(q)\int_{{\mathbb{D}}}{(1-\left\vert{sz}\right\vert ^{2})^{2q-1}\log
 ^{+}\left\vert{f(sz)}\right\vert }.$\ \par 
By proposition~\ref{CI4}, with $\displaystyle p=0,$ we get:\ \par 
\quad \quad \quad $\displaystyle -\int_{{\mathbb{D}}}{\Delta g_{C,s}(z)\log ^{-}\left\vert{f(sz)}\right\vert
 }\leq 4\lbrack 1+2sq\rbrack \int_{{\mathbb{D}}}{(1-\left\vert{sz}\right\vert
 ^{2})^{2q}\log ^{+}\left\vert{f(sz)}\right\vert }.$\ \par 
So, adding:\ \par 
\quad \quad \quad $\displaystyle \ \int_{{\mathbb{D}}}{\Delta g_{C,s}(z)\log \left\vert{f(sz)}\right\vert
 }\leq c(q)\int_{{\mathbb{D}}}{(1-\left\vert{sz}\right\vert ^{2})^{2q-1}\log
 ^{+}\left\vert{f(sz)}\right\vert }.$\ \par 
\quad By proposition~\ref{eOP4} with $\displaystyle p=0,$ we get:\ \par 
\quad \quad \quad $\displaystyle \triangle g_{A,s}(z)\lesssim d(sz,E)^{2q-1},$\ \par 
hence\ \par 
\quad \quad \quad $\displaystyle \ \int_{{\mathbb{D}}}{\Delta g_{A,s}(z)\log ^{+}\left\vert{f(sz)}\right\vert
 }\lesssim \int_{{\mathbb{D}}}{d(sz,E)^{2q-1}\log ^{+}\left\vert{f(sz)}\right\vert
 }.$\ \par 
\quad By proposition~\ref{eOP7} we get:\ \par 
\quad \quad \quad $\displaystyle -\int_{{\mathbb{D}}}{\triangle g_{A,s}(z)\log
 ^{-}\left\vert{f(sz)}\right\vert \lesssim }\int_{{\mathbb{D}}}{(1-\left\vert{sz}\right\vert
 ^{2})^{2q-1}\log ^{+}\left\vert{f(sz)}\right\vert }+$\ \par 
\quad \quad \quad \quad \quad \quad \quad \quad \quad $\displaystyle +2{\times}4^{q}(1-u^{2})^{-2q}\int_{D(0,u)}{(1-\left\vert{z}\right\vert
 ^{2})^{-1/2}d(sz,E)^{2q}\log ^{+}\left\vert{fsz}\right\vert }+$\ \par 
\quad \quad \quad \quad \quad \quad \quad \quad \quad $\displaystyle +2(1-u^{2})^{1/4}u^{-2}\int_{{\mathbb{D}}\backslash
 D(0,u)}{(1-\left\vert{z}\right\vert ^{2})^{-3/4}\left\vert{z}\right\vert
 ^{2}\varphi _{A}(sz)\log ^{-}\left\vert{fsz}\right\vert }.$\ \par 
\quad So adding we get\ \par 
\quad \quad \quad $\displaystyle \ \int_{{\mathbb{D}}}{\Delta g_{A,s}(z)\log \left\vert{f(sz)}\right\vert
 }\lesssim \int_{{\mathbb{D}}}{d(sz,E)^{2q-1}\log ^{+}\left\vert{f(sz)}\right\vert
 }+\int_{{\mathbb{D}}}{(1-\left\vert{sz}\right\vert ^{2})^{2q-1}\log
 ^{+}\left\vert{f(sz)}\right\vert }+$\ \par 
\quad \quad \quad \quad \quad \quad \quad \quad \quad $\displaystyle +2{\times}4^{q}(1-u^{2})^{-2q}\int_{D(0,u)}{(1-\left\vert{z}\right\vert
 ^{2})^{-1/2}d(sz,E)^{2q}\log ^{+}\left\vert{fsz}\right\vert }+$\ \par 
\quad \quad \quad \quad \quad \quad \quad \quad \quad $\displaystyle +2(1-u^{2})^{1/4}u^{-2}\int_{{\mathbb{D}}\backslash
 D(0,u)}{(1-\left\vert{z}\right\vert ^{2})^{-3/4}\left\vert{z}\right\vert
 ^{2}\varphi _{A}(sz)\log ^{-}\left\vert{fsz}\right\vert }.$\ \par 
\quad Combining these results, we proved:\ \par 

\begin{Prps}
~\label{eP4}We have:\par 
\quad \quad \quad $\displaystyle \ \int_{{\mathbb{D}}}{\Delta g_{s}(z)\log \left\vert{f(sz)}\right\vert
 }\lesssim \int_{{\mathbb{D}}}{d(sz,E)^{2q-1}\log ^{+}\left\vert{f(sz)}\right\vert
 }+\int_{{\mathbb{D}}}{(1-\left\vert{sz}\right\vert ^{2})^{2q-1}\log
 ^{+}\left\vert{f(sz)}\right\vert }+$\par 
\quad \quad \quad \quad \quad \quad \quad \quad \quad $\displaystyle +2{\times}4^{q}(1-u^{2})^{-2q}\int_{D(0,u)}{(1-\left\vert{z}\right\vert
 ^{2})^{-1/2}d(sz,E)^{2q}\log ^{+}\left\vert{fsz}\right\vert }+$\par 
\quad \quad \quad \quad \quad \quad \quad \quad \quad $\displaystyle +4(1-u^{2})^{1/2}\sup _{su<t<s}\int_{{\mathbb{T}}}{\varphi
 _{A}(te^{i\theta })\log ^{-}\left\vert{f(te^{i\theta })d\theta }\right\vert }.$
\end{Prps}
\quad Proof.\ \par 
It remains to deal with the term in $\displaystyle \log ^{-}\left\vert{fsz}\right\vert
 .$ We have, passing in polar coordinates,\ \par 
\quad \quad \quad $\displaystyle \ \int_{{\mathbb{D}}\backslash D(0,u)}{(1-\left\vert{z}\right\vert
 ^{2})^{-3/4}\left\vert{z}\right\vert ^{2}\varphi _{A}(sz)\log
 ^{-}\left\vert{fsz}\right\vert }=$\ \par 
\quad \quad \quad \quad \quad \quad \quad $\displaystyle \ \int_{u}^{1}{(1-\rho ^{2})^{-3/4}\lbrace \int_{{\mathbb{T}}}{\varphi
 _{A}(s\rho e^{i\theta })\log ^{-}\left\vert{fs\rho e^{i\theta
 }}\right\vert }\rho d\rho }\leq $\ \par 
\quad \quad \quad \quad \quad \quad \quad \quad \quad $\displaystyle \leq \sup _{su<t<s}\int_{{\mathbb{T}}}{\varphi
 _{A}(te^{i\theta })\log ^{-}\left\vert{f(te^{i\theta })d\theta
 }\right\vert }\int_{u}^{1}{(1-\rho ^{2})^{-3/4}\rho d\rho }\leq $\ \par 
\quad \quad \quad \quad \quad \quad \quad \quad \quad $\displaystyle \leq 2(1-u^{2})^{1/4}\sup _{su<t<s}\int_{{\mathbb{T}}}{\varphi
 _{A}(te^{i\theta })\log ^{-}\left\vert{f(te^{i\theta })d\theta
 }\right\vert }.$\ \par 
Hence we get\ \par 
\quad \quad \quad $\displaystyle 2(1-u^{2})^{1/4}u^{-2}\int_{{\mathbb{D}}\backslash
 D(0,u)}{(1-\left\vert{z}\right\vert ^{2})^{-3/4}\left\vert{z}\right\vert
 ^{2}\varphi _{A}(sz)\log ^{-}\left\vert{fsz}\right\vert }\leq $\ \par 
\quad \quad \quad \quad \quad \quad \quad $\displaystyle \leq 4(1-u^{2})^{1/2}\sup _{su<t<s}\int_{{\mathbb{T}}}{\varphi
 _{A}(te^{i\theta })\log ^{-}\left\vert{f(te^{i\theta })d\theta
 }\right\vert }$\ \par 
which ends the proof. $\hfill\blacksquare $\ \par 
\ \par 
\quad Now we shall use the relation~(\ref{eP3}) which says:\ \par 
\quad \quad \quad $\displaystyle \ \sum_{a\in Z(f_{s})}{g_{s}(a)}=\int_{{\mathbb{D}}}{\log
 \left\vert{f(sz)}\right\vert \triangle g_{s}(z)}+2\int_{{\mathbb{T}}}{\varphi
 (e^{i\theta })\log \left\vert{f(se^{i\theta })}\right\vert }.$\ \par 
So we have, using proposition~\ref{eP4},\ \par 
\quad \quad \quad $\displaystyle \ \sum_{a\in Z(f_{s})}{g_{s}(a)}\lesssim \int_{{\mathbb{D}}}{d(sz,E)^{2q-1}\log
 ^{+}\left\vert{f(sz)}\right\vert }+\int_{{\mathbb{D}}}{(1-\left\vert{sz}\right\vert
 ^{2})^{2q-1}\log ^{+}\left\vert{f(sz)}\right\vert }+$\ \par 
\quad \quad \quad \quad \quad \quad \quad \quad $\displaystyle +2{\times}4^{q}(1-u^{2})^{-2q}\int_{D(0,u)}{(1-\left\vert{z}\right\vert
 ^{2})^{-1/2}d(sz,E)^{2q}\log ^{+}\left\vert{fsz}\right\vert }+$\ \par 
\quad \quad \quad \quad \quad \quad \quad \quad \quad $\displaystyle +4(1-u^{2})^{1/2}\sup _{su<t<s}\int_{{\mathbb{T}}}{\varphi
 _{A}(te^{i\theta })\log ^{-}\left\vert{f(te^{i\theta })d\theta
 }\right\vert }+$\ \par 
\quad \quad \quad \quad \quad \quad \quad \quad \quad $\displaystyle +2\int_{{\mathbb{T}}}{\varphi (se^{i\theta })\log
 ^{+}\left\vert{f(se^{i\theta })}\right\vert }-2\int_{{\mathbb{T}}}{\varphi
 (se^{i\theta })\log ^{-}\left\vert{f(se^{i\theta })}\right\vert }.$\ \par 
\quad So the "good" term is now $\displaystyle -2\int_{{\mathbb{T}}}{\varphi
 (se^{i\theta })\log ^{-}\left\vert{f(se^{i\theta })}\right\vert }.$\ \par 
We shall set\ \par 
\quad \quad \quad $\displaystyle P_{{\mathbb{T}},+}(t_{0}):=\sup _{0\leq s\leq
 t_{0}}\int_{{\mathbb{T}}}{\varphi (se^{i\theta })\log ^{+}\left\vert{f(se^{i\theta
 })}\right\vert }.$\ \par 
and\ \par 
\quad \quad \quad $\displaystyle P_{{\mathbb{T}},-}(t_{0}):=\sup _{0\leq s\leq
 t_{0}}\int_{{\mathbb{T}}}{\varphi (se^{i\theta })\log ^{-}\left\vert{f(se^{i\theta
 })}\right\vert }.$\ \par 
Because $\displaystyle \gamma (s):=\int_{{\mathbb{T}}}{\varphi
 (se^{i\theta })\log ^{-}\left\vert{f(se^{i\theta })}\right\vert
 }$ is continuous for $\displaystyle s\in \lbrack 0,t_{0}\rbrack
 $ by lemma~\ref{2_CF14} in the appendix,  the $\displaystyle
 \sup $ is achieved for a $\displaystyle s_{0}\in \lbrack 0,t_{0}\rbrack
 $ and we have\ \par 
\quad \quad \quad $\displaystyle P_{{\mathbb{T}},-}(t_{0})=\int_{{\mathbb{T}}}{\varphi
 (s_{0}e^{i\theta })\log ^{-}\left\vert{f(s_{0}e^{i\theta })}\right\vert
 }.$\ \par 
\ \par 
\quad Fix $\displaystyle t_{0}<1$ and set:\ \par 
\quad \quad \quad $\displaystyle P_{{\mathbb{D}},+}(s):=\int_{{\mathbb{D}}}{d(sz,E)^{2q-1}\log
 ^{+}\left\vert{f(sz)}\right\vert }+\int_{{\mathbb{D}}}{(1-\left\vert{sz}\right\vert
 ^{2})^{2q-1}\log ^{+}\left\vert{f(sz)}\right\vert }+$\ \par 
\quad \quad \quad \quad \quad \quad \quad \quad $\displaystyle +2{\times}4^{q}(1-u^{2})^{-2q}\int_{D(0,u)}{(1-\left\vert{z}\right\vert
 ^{2})^{-1/2}d(sz,E)^{2q}\log ^{+}\left\vert{fsz}\right\vert }+$\ \par 
\quad \quad \quad \quad \quad \quad \quad \quad \quad $\displaystyle +2\int_{{\mathbb{T}}}{\varphi (se^{i\theta })\log
 ^{+}\left\vert{f(se^{i\theta })}\right\vert }.$\ \par 
Then we get, with the $\displaystyle s_{0}\leq t_{0}$ associated
 to $\displaystyle t_{0},$\ \par 
\quad \quad \quad $\displaystyle \ \sum_{a\in Z(f_{s})}{g_{t_{0}}(a)}+\sum_{a\in
 Z(f_{s})}{g_{s_{0}}(a)}\lesssim P_{{\mathbb{D}},+}(t_{0})+P_{{\mathbb{D}},+}(s_{0})+$\
 \par 
\quad \quad \quad \quad \quad \quad \quad \quad \quad $\displaystyle +4(1-u^{2})^{1/2}\sup _{t_{0}u<t<t_{0}}\int_{{\mathbb{T}}}{\varphi
 _{A}(te^{i\theta })\log ^{-}\left\vert{f(te^{i\theta })d\theta
 }\right\vert }+$\ \par 
\quad \quad \quad \quad \quad \quad \quad \quad \quad $\displaystyle +4(1-u^{2})^{1/2}\sup _{s_{0}u<t<s_{0}}\int_{{\mathbb{T}}}{\varphi
 _{A}(te^{i\theta })\log ^{-}\left\vert{f(te^{i\theta })d\theta
 }\right\vert }-$\ \par 
\quad \quad \quad \quad \quad \quad \quad \quad \quad $\displaystyle -2\int_{{\mathbb{T}}}{\varphi (t_{0}e^{i\theta
 })\log ^{-}\left\vert{f(t_{0}e^{i\theta })}\right\vert }-2\int_{{\mathbb{T}}}{\varphi
 (s_{0}e^{i\theta })\log ^{-}\left\vert{f(t_{0}e^{i\theta })}\right\vert
 }.$\ \par 
But, because $\displaystyle P_{{\mathbb{T}},-}(t_{0})=\int_{{\mathbb{T}}}{\varphi
 (s_{0}e^{i\theta })\log ^{-}\left\vert{f(s_{0}e^{i\theta })}\right\vert
 }$ and $\displaystyle s_{0}\leq t_{0},$ we get:\ \par 
\quad \quad \quad $\displaystyle \sup _{t_{0}u<t<t_{0}}\int_{{\mathbb{T}}}{\varphi
 _{A}(te^{i\theta })\log ^{-}\left\vert{f(te^{i\theta })d\theta
 }\right\vert }\leq P_{{\mathbb{T}},-}(t_{0})$\ \par 
and\ \par 
\quad \quad \quad $\displaystyle \sup _{s_{0}u<t<s_{0}}\int_{{\mathbb{T}}}{\varphi
 _{A}(te^{i\theta })\log ^{-}\left\vert{f(te^{i\theta })d\theta
 }\right\vert }\leq P_{{\mathbb{T}},-}(t_{0}),$\ \par 
so\ \par 
\quad \quad \quad $\displaystyle 8(1-u^{2})^{1/2}P_{{\mathbb{T}},-}(t_{0})-2\int_{{\mathbb{T}}}{\varphi
 (s_{0}e^{i\theta })\log ^{-}\left\vert{f(t_{0}e^{i\theta })}\right\vert
 }\leq $\ \par 
\quad \quad \quad \quad \quad \quad \quad $\displaystyle \leq (8(1-u^{2})^{1/2}-2)\int_{{\mathbb{T}}}{\varphi
 (s_{0}e^{i\theta })\log ^{-}\left\vert{f(t_{0}e^{i\theta })}\right\vert
 }.$\ \par 
So choosing $\displaystyle u<1$ such that $\displaystyle 8(1-u^{2})^{1/2}-2\leq
 0,$ i.e. $\displaystyle u\geq {\sqrt{\frac{15}{16}}}$ which
 is {\sl independent} of $\displaystyle t_{0},$ we get\ \par 
\quad \begin{equation}  \forall t_{0}<1,\ \sum_{a\in Z(f_{s})}{g_{t_{0}}(a)}\leq
 \sum_{a\in Z(f_{s})}{g_{t_{0}}(a)}+\sum_{a\in Z(f_{s})}{g_{s_{0}}(a)}\lesssim
 P_{{\mathbb{D}},+}(t_{0})+P_{{\mathbb{D}},+}(s_{0}).\label{eP5}\end{equation}\
 \par 
In fact we have, for $\displaystyle u={\sqrt{\frac{15}{16}}},$\ \par 
\quad \quad \quad $\displaystyle 2{\times}4^{q}(1-u^{2})^{-2q}\int_{D(0,u)}{(1-\left\vert{z}\right\vert
 ^{2})^{-1/2}d(sz,E)^{2q}\log ^{+}\left\vert{fsz}\right\vert }\leq $\ \par 
\quad \quad \quad \quad \quad \quad \quad $\leq 2{\times}4^{q}(\frac{16}{15})^{q}\lbrack \int_{{\mathbb{D}}}{d(sz,E)^{2q-1}\log
 ^{+}\left\vert{f(sz)}\right\vert }+\int_{{\mathbb{D}}}{(1-\left\vert{sz}\right\vert
 ^{2})^{2q-1}\log ^{+}\left\vert{f(sz)}\right\vert }\rbrack .$\ \par 
So\ \par 
\quad \quad \quad $\displaystyle P_{{\mathbb{D}},+}(s)\leq c(q)\lbrack \int_{{\mathbb{D}}}{d(sz,E)^{2q-1}\log
 ^{+}\left\vert{f(sz)}\right\vert }+\int_{{\mathbb{D}}}{(1-\left\vert{sz}\right\vert
 ^{2})^{2q-1}\log ^{+}\left\vert{f(sz)}\right\vert }\rbrack +P_{{\mathbb{T}},+}(s).$\
 \par 
\quad So we are lead to the definition, replacing $\displaystyle 2q$ by $q:$\ \par 

\begin{Dfnt}
Let $\displaystyle E=\bar E\subset {\mathbb{T}}.$ We say that
 an holomorphic function $f$ is in the generalised Nevanlinna
 class $\displaystyle {\mathcal{N}}_{d(\cdot ,E)^{q},0}$ if 
 $\displaystyle \exists \delta >0,\ \delta <1$ such that\par 
\quad \quad \quad $\displaystyle \ {\left\Vert{f}\right\Vert}_{{\mathcal{N}}_{d(\cdot
 ,E)^{q},0}}:=\sup _{1-\delta <s<1}\lbrace \int_{{\mathbb{T}}}{d(se^{i\theta
 },E)^{q}\log ^{+}\left\vert{f(se^{i\theta })}\right\vert }+$\par 
\quad \quad \quad \quad \quad \quad \quad $\displaystyle \ +\int_{{\mathbb{D}}}{d(sz,E)^{q-1}\log ^{+}\left\vert{f(sz)}\right\vert
 }+\int_{{\mathbb{D}}}{(1-\left\vert{sz}\right\vert ^{2})^{q-1}\log
 ^{+}\left\vert{f(sz)}\right\vert }\rbrace <\infty .$
\end{Dfnt}
\quad And we proved the Blaschke type condition:\ \par 

\begin{Thrm}
Let $\displaystyle E=\bar E\subset {\mathbb{T}}.$ Suppose $\displaystyle
 q>0$ and $\displaystyle f\in {\mathcal{N}}_{\varphi ,0}({\mathbb{D}})$
 with $\displaystyle \ \left\vert{f(0)}\right\vert =1,$ then\par 
\quad \quad \quad $\displaystyle \ \sum_{a\in Z(f)}{(1-\left\vert{a}\right\vert
 ^{2})\varphi (a)}\leq c(\varphi ){\left\Vert{f}\right\Vert}_{{\mathcal{N}}_{\varphi
 ,0}}.$
\end{Thrm}

\begin{Crll}
~\label{0_NI0}Let $\displaystyle E=\bar E\subset {\mathbb{T}}.$
 Suppose $\displaystyle q>0$ and $\displaystyle f\in {\mathcal{N}}_{d(\cdot
 ,E)^{q},0}({\mathbb{D}})$ with $\displaystyle \ \left\vert{f(0)}\right\vert
 =1,$ then\par 
\quad \quad \quad $\displaystyle \ \sum_{a\in Z(f)}{(1-\left\vert{a}\right\vert
 ^{2})d(a,E)^{q}}\leq c(E,q){\left\Vert{f}\right\Vert}_{{\mathcal{N}}_{d(\cdot
 ,E)^{q},0}}.$
\end{Crll}
\quad Proof.\ \par 
For the theorem we apply inequality~(\ref{eP5}) \ \par 
\quad \quad \quad $\displaystyle \forall t<1,\ \sum_{a\in Z(f_{t})}{g_{t}(a)}\lesssim
 P_{{\mathbb{D}},+}(t)+P_{{\mathbb{D}},+}(s),$\ \par 
and for the corollary we recall that $\displaystyle \varphi (z)\simeq
 d(z,E)^{q}.$ $\hfill\blacksquare $\ \par 

\section{Application : $L^{\infty }$ bounds.}
We shall examine two cases.\ \par 
\ \par 
\quad $\bullet $ {\bf Case }$\displaystyle p>0.$\ \par 
\ \par 
\quad Let $\displaystyle E=\bar E\subset {\mathbb{T}}\ ;$ its Ahern-Clark
 type $\alpha (E)$ is defined the following way:\ \par 
\quad \quad \quad $\displaystyle \alpha (E):=\sup \lbrace \alpha \in {\mathbb{R}}::\left\vert{\lbrace
 t\in {\mathbb{T}}::d(t,E)<x\rbrace }\right\vert ={\mathcal{O}}(x^{\alpha
 }),\ x\rightarrow +0\rbrace ,$\ \par 
where $\displaystyle \ \left\vert{A}\right\vert $ denotes the
 Lebesgue measure of the set $\displaystyle A.$\ \par 
\quad Our hypothesis is\ \par 
\quad \quad \quad $\displaystyle \log \left\vert{f(z)}\right\vert \leq \frac{K}{(1-\left\vert{z}\right\vert
 )^{p}}\frac{1}{d(z,E)^{q}},\ z\in {\mathbb{D}},\ p,q\geq 0.$\ \par 
We want to apply corollary~\ref{eP2} so we have, with $\varphi
 (z):=d(z,E)^{q-\alpha (E)+\epsilon }:$\ \par 
\quad \quad \quad $\displaystyle \ \sum_{a\in Z(f)}{(1-\left\vert{a}\right\vert
 ^{2})^{1+p}\varphi (a)}\leq c(\varphi ){\left\Vert{f}\right\Vert}_{{\mathcal{N}}_{\varphi
 ,p}},$\ \par 
and we shall compute $\displaystyle \ {\left\Vert{f}\right\Vert}_{{\mathcal{N}}_{\varphi
 ,p}},$ i.e.\ \par 
\quad \quad \quad $\displaystyle \ {\left\Vert{f}\right\Vert}_{{\mathcal{N}}_{\varphi
 ,p}}:=\sup _{1-\delta \leq s<1}\int_{{\mathbb{D}}}{(1-\left\vert{z}\right\vert
 ^{2})^{p-1}d(sz,E)^{q-\alpha (E)+\epsilon }\log ^{+}\left\vert{f(sz)}\right\vert
 }.$\ \par 
The hypothesis gives\ \par 
\quad \quad \quad $\displaystyle \forall z\in {\mathbb{D}},\ \log ^{+}\left\vert{f(z)}\right\vert
 \leq \frac{K}{(1-\left\vert{z}\right\vert ^{2})^{p}}\frac{1}{d(z,E)^{q}(z)},$\
 \par 
so we have\ \par 
\quad \quad \quad $\displaystyle \ \int_{{\mathbb{D}}}{(1-\left\vert{z}\right\vert
 ^{2})^{p-1}d(sz,E)^{q-\alpha (E)+\epsilon }\log ^{+}\left\vert{f(sz)}\right\vert
 }\leq K\int_{{\mathbb{D}}}{(1-\left\vert{z}\right\vert ^{2})^{\epsilon
 -1}d(sz,E)^{-\alpha (E)}}.$\ \par 
We set $\displaystyle \Gamma _{n}:=E_{n}{\times}(1-2^{-n},\ 1)$
 and $\gamma _{n}:=\Gamma _{n}\backslash \Gamma _{n+1}.$ Then we get\ \par 
\quad \quad \quad $\displaystyle \ \int_{{\mathbb{D}}}{(1-\left\vert{z}\right\vert
 ^{2})^{\epsilon -1}d(sz,E)^{-\alpha (E)}}=\sum_{n\in {\mathbb{N}}}{\int_{\gamma
 _{n}}{(1-\left\vert{z}\right\vert ^{2})^{\epsilon -1}d(sz,E)^{-\alpha
 (E)}}}\leq $\ \par 
\quad \quad \quad \quad \quad $\displaystyle \ \sum_{n\in {\mathbb{N}}}{2^{-(\epsilon -1)n}2^{n\alpha
 (E)}\int_{\gamma _{n}}{dm(z)}}\leq \sum_{n\in {\mathbb{N}}}{2^{-(\epsilon
 -1)n}2^{n\alpha (E)}\left\vert{E_{n}}\right\vert 2^{-n}}=\sum_{n\in
 {\mathbb{N}}}{2^{-\epsilon n}}=:c(\epsilon )<\infty $\ \par 
because $\epsilon >0.$\ \par 
\quad So corollary~\ref{eP2} gives\ \par 
\quad \quad \quad $\displaystyle \ \sum_{a\in Z(f)}{(1-\left\vert{a}\right\vert
 ^{2})^{1+p}d(a,E)^{q-\alpha (E)+\epsilon )}}\leq c(p,q,\epsilon
 ){\left\Vert{f}\right\Vert}_{{\mathcal{N}}_{\varphi ,p}},$\ \par 
hence we get\ \par 
\quad \quad \quad $\displaystyle \ \sum_{a\in Z(f)}{(1-\left\vert{a}\right\vert
 ^{2})^{1+p}d(a,E)^{q-\alpha (E)+\epsilon }}\leq c(p,q,\epsilon
 ){\left\Vert{f}\right\Vert}_{{\mathcal{N}}_{\varphi ,p}}\leq
 Kc(\epsilon )c(p,q,\epsilon ).$\ \par 
So we proved:\ \par 

\begin{Thrm}
Suppose that $\displaystyle f\in {\mathcal{H}}({\mathbb{D}}),\
 \left\vert{f(0)}\right\vert =1$ and\par 
\quad \quad \quad $\displaystyle \forall z\in {\mathbb{D}},\ \log ^{+}\left\vert{f(z)}\right\vert
 \leq \frac{K}{(1-\left\vert{z}\right\vert ^{2})^{p}}\frac{1}{d(z,E)^{q}},$\par 
then we have, with any $\displaystyle \epsilon >0,$\par 
\quad \quad \quad $\displaystyle \ \sum_{a\in Z(f)}{(1-\left\vert{a}\right\vert
 ^{2})^{1+p}d(a,E)^{(q-\alpha (E)+\epsilon )_{+}}}\leq c(p,q,\epsilon )K.$
\end{Thrm}
\ \par 
{\large       $\bullet $} {\bf Case} $\displaystyle p=0.$\ \par 
\ \par 
\quad For this case we want to apply corollary~\ref{0_NI0}. So let\ \par 
\quad \quad \quad $\displaystyle \forall z\in {\mathbb{D}},\ \log ^{+}\left\vert{f(z)}\right\vert
 \leq K\frac{1}{d(z,E)^{q}}.$\ \par 
\quad We have\ \par 

\begin{Thrm}
Suppose that $\displaystyle f\in {\mathcal{H}}({\mathbb{D}}),\
 \left\vert{f(0)}\right\vert =1$ and\par 
\quad \quad \quad $\displaystyle \forall z\in {\mathbb{D}},\ \log ^{+}\left\vert{f(z)}\right\vert
 \leq K\frac{1}{d(z,E)^{q}},$\par 
then\par 
\quad \quad \quad $\displaystyle \ \sum_{a\in Z(f)}{(1-\left\vert{a}\right\vert
 ^{2})d(a,E)^{(q-\alpha (E)+\epsilon )_{+}}}\leq c(q,\epsilon )K.$
\end{Thrm}
\quad Proof.\ \par 
We have to verify\ \par 
\quad \quad \quad $\displaystyle \sup _{1-\delta \leq s<1}\int_{{\mathbb{T}}}{d(se^{i\theta
 },E)^{q-\alpha (E)+\epsilon }\log ^{+}\left\vert{f(se^{i\theta
 })}\right\vert d\theta }<\infty $\ \par 
and\ \par 
\quad \quad \quad $\displaystyle \sup _{1-\delta \leq s<1}\int_{{\mathbb{D}}}{d(sz,E)^{q-\alpha
 (E)-1+\epsilon }\log ^{+}\left\vert{f(sz)}\right\vert }<\infty .$\ \par 
\quad For the first one, we have\ \par 
\quad \quad \quad $\displaystyle \ \int_{{\mathbb{T}}}{d(se^{i\theta },E)^{q-\alpha
 (E)+\epsilon }\log ^{+}\left\vert{f(se^{i\theta })}\right\vert
 }\leq K\int_{{\mathbb{T}}}{d(se^{i\theta },E)^{\epsilon -\alpha (E)}}.$\ \par 
Set $\displaystyle E_{n}:=\lbrace x\in {\mathbb{T}}::d(x,E)\geq
 2^{-n}\rbrace $ and $\displaystyle F_{n}:=E_{n}\backslash E_{n+1}.$
 we have\ \par 
\quad \quad \quad $\displaystyle \ \int_{{\mathbb{T}}}{d(se^{i\theta },E)^{\epsilon
 -\alpha (E)}d\theta }=\ \sum_{n\in {\mathbb{N}}}{\int_{F_{n}}{d(se^{i\theta
 },E)^{\epsilon -\alpha (E)}d\theta }}\leq \sum_{n\in {\mathbb{N}}}{\int_{F_{n}}{2^{-n(\epsilon
 -\alpha (E))}d\theta }}\leq $\ \par 
\quad \quad \quad \quad \quad \quad \quad \quad \quad $\displaystyle \ \sum_{n\in {\mathbb{N}}}{2^{-n(\epsilon -\alpha
 (E))}\int_{F_{n}}{d\theta }}\leq \sum_{n\in {\mathbb{N}}}{2^{-n\epsilon
 }}<\infty $\ \par 
by the very definition of $\displaystyle \alpha (E)$ and because
 $\displaystyle \epsilon >0.$\ \par 
\quad For the second one we set $\displaystyle \Gamma _{n}:=E_{n}{\times}(1-2^{-n},\
 1)$ and $\gamma _{n}:=\Gamma _{n}\backslash \Gamma _{n+1}.$\ \par 
We get\ \par 
\quad \quad \quad $\displaystyle \ \int_{{\mathbb{D}}}{d(sz,E)^{q-\alpha (E)-1+\epsilon
 }\log ^{+}\left\vert{f(sz)}\right\vert }=\int_{{\mathbb{D}}}{d(sz,E)^{q-\alpha
 (E)-1+\epsilon }\log ^{+}\left\vert{f(sz)}\right\vert }\leq
 \int_{{\mathbb{D}}}{d(sz,E)^{\epsilon -\alpha (E)-1}}.$\ \par 
But\ \par 
\quad \quad $\displaystyle \ \sum_{n\in {\mathbb{N}}}{\int_{\gamma _{n}}{d(sz,E)^{\epsilon
 -\alpha (E)-1}}}\leq \sum_{n\in {\mathbb{N}}}{\int_{\gamma _{n}}{2^{-n(\epsilon
 -\alpha (E)-1)}}}\leq \sum_{n\in {\mathbb{N}}}{2^{-n(\epsilon
 -\alpha (E)-1)}\int_{\gamma _{n}}{dm(z)}}\leq $\ \par 
\quad \quad \quad \quad \quad \quad \quad \quad \quad $\displaystyle \ \sum_{n\in {\mathbb{N}}}{2^{-n(\epsilon -\alpha
 (E)-1)}\left\vert{E_{n}}\right\vert {\times}(2^{-n})}\leq \sum_{n\in
 {\mathbb{N}}}{2^{-n\epsilon }}<\infty ,$\ \par 
because $\displaystyle \epsilon >0.$ We end the proof as in the
 case $\displaystyle p>0.$ $\hfill\blacksquare $\ \par 
\quad These results give alternative proofs of some of the results
 by Favorov \&  Golinskii~\cite{Golinskii12}.\ \par 

\section{Appendix.}
\quad When there is no ambiguities, we shall forget the index $\displaystyle
 j.$\ \par 

\begin{Lmm}
~\label{CI23}We have\par 
\quad \quad \quad $\displaystyle \bar \partial (\psi _{j})^{q}(z)=q\frac{(z-\alpha
 _{j})\left\vert{z-\alpha _{j}}\right\vert ^{2q-2}\left\vert{z-\beta
 _{j}}\right\vert ^{2q}}{\delta _{j}^{2q}}+q\frac{(z-\beta _{j})\left\vert{z-\alpha
 _{j}}\right\vert ^{2q}\left\vert{z-\beta _{j}}\right\vert ^{2q-2}}{\delta
 _{j}^{2q}}.$\par 
And\par 
\quad \quad \quad $\displaystyle \partial \bar \partial (\psi _{j})^{q}=q^{2}\frac{\left\vert{z-\alpha
 _{j}}\right\vert ^{2q-2}\left\vert{z-\beta _{j}}\right\vert
 ^{2q-2}}{\delta _{j}^{2q}}\lbrace \left\vert{z-\alpha _{j}}\right\vert
 ^{2}+\left\vert{z-\beta _{j}}\right\vert ^{2}+2\Re \lbrack (z-\alpha
 _{j})(\bar z-\bar \beta _{j})\rbrack \rbrace .$
\end{Lmm}
\quad Proof.\ \par 
We have $\displaystyle \forall z\in \Gamma ,\ \psi (z)^{q}=\frac{\left\vert{z-\alpha
 }\right\vert ^{2q}\left\vert{z-\beta }\right\vert ^{2q}}{\delta
 ^{2q}}$ so\ \par 
\quad \quad \quad $\displaystyle \bar \partial (\psi )^{q}(z)=q\frac{(z-\alpha
 )\left\vert{z-\alpha }\right\vert ^{2q-2}\left\vert{z-\beta
 }\right\vert ^{2q}}{\delta ^{2q}}+q\frac{(z-\beta )\left\vert{z-\alpha
 }\right\vert ^{2q}\left\vert{z-\beta }\right\vert ^{2q-2}}{\delta
 ^{2q}}.$\ \par 
And\ \par 
\quad \quad \quad $\displaystyle \partial \bar \partial (\psi )^{q}(z)=q^{2}\frac{\left\vert{z-\alpha
 }\right\vert ^{2q-2}\left\vert{z-\beta }\right\vert ^{2q}}{\delta
 ^{2q}}+q^{2}\frac{\left\vert{z-\alpha }\right\vert ^{2q}\left\vert{z-\beta
 }\right\vert ^{2q-2}}{\delta ^{2q}}+$\ \par 
\quad \quad \quad \quad \quad \quad \quad \quad \quad $\displaystyle +2q^{2}\frac{\left\vert{z-\alpha }\right\vert
 ^{2q-2}\left\vert{z-\beta }\right\vert ^{2q-2}}{\delta ^{2q}}\Re
 \lbrack (z-\alpha )(\bar z-\bar \beta )\rbrack =$\ \par 
\quad \quad \quad \quad \quad \quad \quad \quad $\displaystyle =q^{2}\frac{\left\vert{z-\alpha }\right\vert ^{2q-2}\left\vert{z-\beta
 }\right\vert ^{2q-2}}{\delta ^{2q}}\lbrace \left\vert{z-\beta
 }\right\vert ^{2}+\left\vert{z-\alpha }\right\vert ^{2}+2\Re
 \lbrack (z-\alpha )(\bar z-\bar \beta )\rbrack \rbrace .$ $\hfill\blacksquare
 $\ \par 

\begin{Rmrq}
~\label{CI24}We notice that: $\displaystyle \partial \bar \partial
 (\psi _{j})^{q}(z)\geq 0$ because\par 
\quad \quad \quad $\displaystyle \ \left\vert{z-\beta }\right\vert ^{2}+\left\vert{z-\alpha
 }\right\vert ^{2}+2\Re \lbrack (z-\alpha )(\bar z-\bar \beta
 )\rbrack \geq \left\vert{z-\beta }\right\vert ^{2}+\left\vert{z-\alpha
 }\right\vert ^{2}-2\left\vert{z-\alpha }\right\vert \left\vert{z-\beta
 }\right\vert \geq 0.$
\end{Rmrq}

\begin{Lmm}
~\label{eOP1}If $\displaystyle \eta _{j}'\neq 0$ or if $\eta
 ''_{j}\neq 0,$ we have:\par 
\quad \quad \quad $\displaystyle \forall z\in \Gamma _{j},\ 2(1-\left\vert{z}\right\vert
 ^{2})^{2}\leq \psi _{j}(z)\leq 3(1-\left\vert{z}\right\vert ^{2})^{2}.$
\end{Lmm}
\quad Proof.\ \par 
If $\eta _{j}'\neq 0$ we have\ \par 
\quad $\bullet $ $\chi '_{j}\neq 0$ which implies $\displaystyle 2\leq
 \frac{\psi _{j}(z)}{(1-\left\vert{z}\right\vert ^{2})^{2}}\leq 3$ hence\ \par 
\quad \quad \quad $\displaystyle 2(1-\left\vert{z}\right\vert ^{2})^{2}\leq \psi
 _{j}(z)\leq 3(1-\left\vert{z}\right\vert ^{2})^{2}.$\ \par 
The same for the second derivatives, which ends the proof of
 the lemma. $\hfill\blacksquare $\ \par 

\begin{Lmm}
~\label{eOP0}We have\par 
\quad \quad \quad $\displaystyle \ \left\vert{\bar \partial \lbrack \chi (\frac{\left\vert{z-\alpha
 }\right\vert ^{2}}{(1-\left\vert{z}\right\vert ^{2})^{2}})\rbrack
 }\right\vert \leq 9\left\vert{\chi '()}\right\vert (1-\left\vert{z}\right\vert
 ^{2})^{-1}.$\par 
and\par 
\quad \quad \quad $\displaystyle \ \left\vert{\partial \bar \partial \chi }\right\vert
 \lesssim (\left\vert{\chi '}\right\vert +\left\vert{\chi ''}\right\vert
 )(1-\left\vert{z}\right\vert ^{2})^{-2}.$
\end{Lmm}
\quad Proof.\ \par 
We have\ \par 
\quad \quad \quad $\displaystyle \bar \partial \lbrack \chi (\frac{\left\vert{z-\alpha
 _{j}}\right\vert ^{2}}{(1-\left\vert{z}\right\vert ^{2})^{2}})\rbrack
 =\chi '()\bar \partial \lbrack \frac{\left\vert{z-\alpha _{j}}\right\vert
 ^{2}}{(1-\left\vert{z}\right\vert ^{2})^{2}}\rbrack =\chi '()\lbrack
 \frac{(z-\alpha _{j})}{(1-\left\vert{z}\right\vert ^{2})^{2}}+2z\frac{\left\vert{z-\alpha
 _{j}}\right\vert ^{2}}{(1-\left\vert{z}\right\vert ^{2})^{3}}\rbrack .$\ \par 
But if $\chi '()\neq 0$ then $\displaystyle 2\leq \frac{\left\vert{z-\alpha
 _{j}}\right\vert ^{2}}{(1-\left\vert{z}\right\vert ^{2})^{2}}\leq
 3$ hence\ \par 
\quad \quad \quad $\displaystyle \ \left\vert{\bar \partial \lbrack \chi (\frac{\left\vert{z-\alpha
 _{j}}\right\vert ^{2}}{(1-\left\vert{z}\right\vert ^{2})^{2}})\rbrack
 }\right\vert \leq \left\vert{\chi '()}\right\vert \lbrack {\sqrt{3}}+6\rbrack
 (1-\left\vert{z}\right\vert ^{2})^{-1}\leq 9\left\vert{\chi
 '()}\right\vert (1-\left\vert{z}\right\vert ^{2})^{-1}.$\ \par 
Now\ \par 
\quad \quad \quad $\displaystyle \partial \bar \partial \lbrack \chi (\frac{\left\vert{z-\alpha
 _{j}}\right\vert ^{2}}{(1-\left\vert{z}\right\vert ^{2})^{2}})\rbrack
 =\partial \lbrace \chi '()\lbrack \frac{(z-\alpha _{j})}{(1-\left\vert{z}\right\vert
 ^{2})^{2}}+2z\frac{\left\vert{z-\alpha _{j}}\right\vert ^{2}}{(1-\left\vert{z}\right\vert
 ^{2})^{3}}\rbrack \rbrace =$\ \par 
\quad \quad \quad \quad \quad \quad \quad $\displaystyle =\partial \lbrace \chi '()\rbrace \lbrack \frac{(z-\alpha
 _{j})}{(1-\left\vert{z}\right\vert ^{2})^{2}}+2z\frac{\left\vert{z-\alpha
 _{j}}\right\vert ^{2}}{(1-\left\vert{z}\right\vert ^{2})^{3}}\rbrack
 +\chi '()\partial \lbrace \frac{(z-\alpha _{j})}{(1-\left\vert{z}\right\vert
 ^{2})^{2}}+2z\frac{\left\vert{z-\alpha _{j}}\right\vert ^{2}}{(1-\left\vert{z}\right\vert
 ^{2})^{3}}\rbrace .$\ \par 
And\ \par 
\quad \quad \quad $\partial \chi '()=\chi "()\lbrack \frac{(\bar z-\bar \alpha
 _{j})}{(1-\left\vert{z}\right\vert ^{2})^{2}}+2\bar z\frac{\left\vert{z-\alpha
 _{j}}\right\vert ^{2}}{(1-\left\vert{z}\right\vert ^{2})^{3}}\rbrack .$\ \par 
so\ \par 
\quad \quad \quad $\displaystyle \ \left\vert{\partial \chi '()}\right\vert \leq
 9\left\vert{\chi ''()}\right\vert (1-\left\vert{z}\right\vert
 ^{2})^{-1}.$\ \par 
And a straightforward computation gives\ \par 
\quad \quad \quad $\displaystyle \ \left\vert{\partial \lbrace \frac{(z-\alpha
 _{j})}{(1-\left\vert{z}\right\vert ^{2})^{2}}+2z\frac{\left\vert{z-\alpha
 _{j}}\right\vert ^{2}}{(1-\left\vert{z}\right\vert ^{2})^{3}}\rbrace
 }\right\vert \lesssim (1-\left\vert{z}\right\vert ^{2})^{-2}.$\ \par 
So the lemma is proved. $\hfill\blacksquare $\ \par 

\begin{Lmm}
~\label{CI6}Let $\displaystyle \eta \in {\mathbb{T}},$ then we
 have $\displaystyle \Re (\bar z(z-\eta ))\leq 0$ iff $\displaystyle
 z\in {\mathbb{D}}\cap D(\frac{\eta }{2},\ \frac{1}{2}).$
\end{Lmm}
\quad Proof.\ \par 
We set $\displaystyle z=\eta t,$ then we have\ \par 
\quad \quad \quad $\displaystyle \bar z(z-\eta )=\bar \eta \bar t(\eta t-\eta )=\bar
 t(t-1).$\ \par 
Hence\ \par 
\quad \quad \quad $\displaystyle \Re (\bar z(z-\eta ))=\Re (\bar t(t-1))=\Re (r^{2}-re^{i\theta
 })=r^{2}-r\cos \theta .$\ \par 
Hence with $\displaystyle t=x+iy=re^{i\theta ,\ }x=r\cos \theta
 ,\ y=r\sin \theta ,$ we get\ \par 
\quad \quad \quad $\displaystyle \Re (\bar t(t-1))\leq 0\iff x^{2}+y^{2}-x\leq 0$\ \par 
which means $\displaystyle (x,y)\in D(\frac{1}{2},\ \frac{1}{2})$
 hence $\displaystyle z\in {\mathbb{D}}\cap D(\frac{\eta }{2},\
 \frac{1}{2}).$ $\hfill\blacksquare $\ \par 

\begin{Lmm}
~\label{CI2} (Substitution 1) We have, for $\delta >0$ and $\displaystyle
 u\in \rbrack 0,1\lbrack ,$ and $\displaystyle \ \left\vert{f(0)}\right\vert
 =1,$\par 
\quad \quad \quad $\displaystyle \ \int_{{\mathbb{D}}}{(1-\left\vert{z}\right\vert
 ^{2})^{p-1+\delta }(1-\left\vert{sz}\right\vert ^{2})^{2q}\log
 ^{-}\left\vert{f(sz)}\right\vert }\leq $\par 
\quad \quad \quad \quad \quad \quad \quad $\displaystyle \leq (1-u^{2})^{\delta }u^{-2}\int_{{\mathbb{D}}}{(1-\left\vert{z}\right\vert
 ^{2})^{p-1}\left\vert{z}\right\vert ^{2}(1-\left\vert{sz}\right\vert
 ^{2})^{2q}\log ^{-}\left\vert{f(sz)}\right\vert }+$\par 
\quad \quad \quad \quad \quad \quad \quad \quad \quad \quad $\displaystyle +\int_{{\mathbb{D}}}{(1-\left\vert{z}\right\vert
 ^{2})^{p-1}(1-\left\vert{sz}\right\vert ^{2})^{2q}\log ^{+}\left\vert{f(sz)}\right\vert
 }.$
\end{Lmm}
\quad Proof.\ \par 
We have\ \par 
\quad \quad \quad $\displaystyle \ \int_{{\mathbb{D}}}{(1-\left\vert{z}\right\vert
 ^{2})^{p-1+\delta }(1-\left\vert{sz}\right\vert ^{2})^{2q}\log
 ^{-}\left\vert{f(sz)}\right\vert }=$\ \par 
\quad \quad \quad \quad \quad \quad \quad $\displaystyle =\int_{D(0,u)}{(1-\left\vert{z}\right\vert ^{2})^{p-1+\delta
 }(1-\left\vert{sz}\right\vert ^{2})^{2q}\log ^{-}\left\vert{f(sz)}\right\vert
 }+$\ \par 
\quad \quad \quad \quad \quad \quad \quad \quad \quad $\displaystyle \ +\int_{{\mathbb{D}}\backslash D(0,u)}{(1-\left\vert{z}\right\vert
 ^{2})^{p-1+\delta }(1-\left\vert{sz}\right\vert ^{2})^{2q}\log
 ^{-}\left\vert{f(sz)}\right\vert }.$\ \par 
\quad For the first term, passing in polar coordinates, we get\ \par 
\quad \quad \quad \begin{equation}  I_{1}=\int_{0}^{u}{(1-\rho ^{2})^{p-1+\delta
 }(1-s^{2}\rho ^{2})^{2q}\lbrace \int_{{\mathbb{T}}}{\log ^{-}\left\vert{f(s\rho
 e^{i\theta })}\right\vert }\rbrace \rho d\rho }.\label{CI1}\end{equation}\
 \par 
The subharmonicity of $\displaystyle \log \left\vert{f(sz)}\right\vert
 $ gives\ \par 
\quad \quad \quad $\displaystyle 0=\log \left\vert{f(0)}\right\vert \leq \int_{{\mathbb{T}}}{\log
 \left\vert{f(s\rho e^{i\theta })}\right\vert }=\int_{{\mathbb{T}}}{\log
 ^{+}\left\vert{f(s\rho e^{i\theta })}\right\vert }-\int_{{\mathbb{T}}}{\log
 ^{-}\left\vert{f(s\rho e^{i\theta })}\right\vert },$\ \par 
hence\ \par 
\quad \quad \quad $\displaystyle \ \int_{{\mathbb{T}}}{\log ^{-}\left\vert{f(s\rho
 e^{i\theta })}\right\vert }\leq \int_{{\mathbb{T}}}{\log ^{+}\left\vert{f(s\rho
 e^{i\theta })}\right\vert }.$\ \par 
Putting it in~(\ref{CI1}) we get\ \par 
\quad \quad \quad $\displaystyle I_{1}\leq \int_{0}^{u}{(1-\rho ^{2})^{p-1+\delta
 }(1-s^{2}\rho ^{2})^{2q}\lbrace \int_{{\mathbb{T}}}{\log ^{+}\left\vert{f(s\rho
 e^{i\theta })}\right\vert }\rbrace \rho d\rho }\leq $\ \par 
\quad \quad \quad \quad \begin{equation}  \leq \int_{D(0,u)}{(1-\left\vert{z}\right\vert
 ^{2})^{p-1+\delta }(1-\left\vert{sz}\right\vert ^{2})^{2q}\log
 ^{+}\left\vert{f(sz)}\right\vert }.\label{0_NI1}\end{equation}\ \par 
\quad For the second term, we have\ \par 
\quad \quad \quad $\displaystyle I_{2}:=\int_{{\mathbb{D}}\backslash D(0,u)}{(1-\left\vert{z}\right\vert
 ^{2})^{p-1+\delta }(1-\left\vert{sz}\right\vert ^{2})^{2q}\log
 ^{-}\left\vert{f(sz)}\right\vert }\leq $\ \par 
\quad \quad \quad \quad \quad \quad \quad $\displaystyle \leq (1-u^{2})^{\delta }u^{-2}\int_{{\mathbb{D}}\backslash
 D(0,u)}{(1-\left\vert{z}\right\vert ^{2})^{p-1}\left\vert{z}\right\vert
 ^{2}(1-\left\vert{sz}\right\vert ^{2})^{2q}\log ^{-}\left\vert{f(sz)}\right\vert
 }.$\ \par 
This ends the proof. $\hfill\blacksquare $\ \par 

\begin{Lmm}
~\label{3_CIZ5}(Substitution 2) We have, for $\displaystyle \delta
 >0$ and any $\displaystyle u,\ 0<u<1,$\par 
\quad \quad \quad $\displaystyle \ \int_{{\mathbb{D}}}{(1-\left\vert{z}\right\vert
 ^{2})^{p-1+\delta }\varphi _{A}(sz)\log ^{-}\left\vert{fsz}\right\vert
 }\leq $\par 
\quad \quad \quad \quad \quad \quad \quad \quad \quad \quad $\displaystyle \leq 4^{q}(1-u^{2})^{-2q}\int_{D(0,u)}{(1-\left\vert{z}\right\vert
 ^{2})^{p-1+\delta }\varphi _{A}(sz)\log ^{+}\left\vert{f(sz)}\right\vert
 }+$\par 
\quad \quad \quad \quad \quad \quad \quad \quad \quad \quad \quad $\displaystyle +(1-u^{2})^{\delta /2}u^{-2}\int_{{\mathbb{D}}\backslash
 D(0,u)}{(1-\left\vert{z}\right\vert ^{2})^{p+\delta /2-1}\left\vert{z}\right\vert
 ^{2}\varphi _{A}(sz)\log ^{-}\left\vert{f(sz)}\right\vert }.$
\end{Lmm}
\quad Proof.\ \par 
We have:\ \par 
\quad \quad \quad $\displaystyle \ \int_{{\mathbb{D}}}{(1-\left\vert{z}\right\vert
 ^{2})^{p-1+\delta }\varphi _{A}(sz)\log ^{-}\left\vert{f(sz)}\right\vert
 }=$\ \par 
\quad \quad \quad \quad \quad \quad \quad $\displaystyle =\int_{D(0,u)}{(1-\left\vert{z}\right\vert ^{2})^{p-1+\delta
 }\varphi _{A}(sz)\log ^{-}\left\vert{f(sz)}\right\vert }+$\ \par 
\quad \quad \quad \quad \quad \quad \quad \quad $\displaystyle +\int_{{\mathbb{D}}\backslash D(0,u)}{(1-\left\vert{z}\right\vert
 ^{2})^{p-1+\delta }\varphi _{A}(sz)\log ^{-}\left\vert{f(sz)}\right\vert
 }.$\ \par 
For the first term we get\ \par 
\quad \quad \quad $\displaystyle \ I_{1}:=\int_{D(0,u)}{(1-\left\vert{z}\right\vert
 ^{2})^{p-1+\delta }\varphi _{A}(sz)\log ^{-}\left\vert{f(sz)}\right\vert
 }.$\ \par 
Because $\displaystyle \forall \alpha \in {\mathbb{T}},\ \left\vert{sz-\alpha
 }\right\vert \leq 2$ we get $\displaystyle \varphi _{A,j}(sz)=\eta
 _{j}(z)\frac{\left\vert{z-\alpha _{j}}\right\vert ^{2q}\left\vert{z-\beta
 _{j}}\right\vert ^{2q}}{\delta _{j}^{2q}},$ hence, in order
 to have $\eta _{j}(z)\neq 0,$ we have $\displaystyle \ \left\vert{z-\alpha
 }\right\vert ^{2}\geq 2(1-\left\vert{z}\right\vert ^{2})^{2}$
 and $\displaystyle \ \left\vert{z-\beta }\right\vert ^{2}\geq
 \lambda (1-\left\vert{z}\right\vert ^{2})^{2}.$ But, with~(\ref{eOPC1}),\ \par 
\quad \quad \quad $\displaystyle \psi _{j}(z):=\frac{\left\vert{z-\alpha _{j}}\right\vert
 ^{2}\left\vert{z-\beta _{j}}\right\vert ^{2}}{\delta _{j}^{2}}\leq
 2\left\vert{z-\alpha _{j}}\right\vert ^{2}\leq 4$\ \par 
hence\ \par 
\quad \quad \quad $\displaystyle \varphi _{A,j}(sz)=\eta _{j}(z)\psi _{j}(z)^{q}\leq
 4^{q}.$\ \par 
So we get\ \par 
\quad \quad \quad $\displaystyle I_{1}\leq 4^{q}\int_{D(0,u)}{(1-\left\vert{z}\right\vert
 ^{2})^{p-1+\delta }\log ^{-}\left\vert{f(sz)}\right\vert }.$\ \par 
and we can apply inequality~(\ref{0_NI1}) to get\ \par 
\quad \quad \quad $\displaystyle I_{1}\leq 4^{q}\int_{D(0,u)}{(1-\left\vert{z}\right\vert
 ^{2})^{p-1+\delta }\log ^{+}\left\vert{f(sz)}\right\vert }\ ;$\ \par 
but\ \par 
\quad \quad \quad $\displaystyle \forall z\in D(0,u),\ 1\leq (1-u^{2})^{-2q}\varphi
 _{A}(sz)$\ \par 
so\ \par 
\quad \quad \quad $\displaystyle I_{1}\leq 4^{q}(1-u^{2})^{-2q}\int_{D(0,u)}{(1-\left\vert{z}\right\vert
 ^{2})^{p-1+\delta }\varphi _{A}(sz)\log ^{+}\left\vert{f(sz)}\right\vert
 }\leq $\ \par 
\quad \quad \quad \quad \quad \quad \quad \quad \quad \quad \quad $\displaystyle \leq 4^{q}(1-u^{2})^{-2q}P_{+}(\delta ,u).$\ \par 
\ \par 
\quad For the second one\ \par 
\quad \quad \quad $\displaystyle \ I_{2}:=\int_{{\mathbb{D}}\backslash D(0,u)}{(1-\left\vert{z}\right\vert
 ^{2})^{p-1+\delta }\varphi _{A}(sz)\log ^{-}\left\vert{f(sz)}\right\vert
 }\leq $\ \par 
\quad \quad \quad \quad \quad \quad \quad \quad \quad $\displaystyle \leq (1-u^{2})^{\delta /2}u^{-2}\int_{{\mathbb{D}}\backslash
 D(0,u)}{(1-\left\vert{z}\right\vert ^{2})^{p+\delta /2-1}\left\vert{z}\right\vert
 ^{2}\varphi _{A}(sz)\log ^{-}\left\vert{f(sz)}\right\vert }.$\ \par 
\quad Adding we get\ \par 
$\displaystyle \ \int_{{\mathbb{D}}}{(1-\left\vert{z}\right\vert
 ^{2})^{p-1+\delta }\varphi _{A}(sz)\log ^{-}\left\vert{fsz}\right\vert
 }\leq 4^{q}(1-u^{2})^{-2q}\int_{D(0,u)}{(1-\left\vert{z}\right\vert
 ^{2})^{p-1+\delta }\log ^{+}\left\vert{f(sz)}\right\vert }+$\ \par 
\quad \quad \quad \quad \quad \quad \quad \quad \quad \quad \quad $\displaystyle +(1-u^{2})^{\delta /2}u^{-2}\int_{{\mathbb{D}}\backslash
 D(0,u)}{(1-\left\vert{z}\right\vert ^{2})^{p+\delta /2-1}\left\vert{z}\right\vert
 ^{2}\varphi _{A}(sz)\log ^{-}\left\vert{f(sz)}\right\vert }.$\ \par 
Which ends the proof of the lemma. $\hfill\blacksquare $\ \par 

\begin{Lmm}
~\label{2_CF14}Let $\varphi $ be a continuous function in the
 unit disc $\displaystyle {\mathbb{D}}.$ We have that:\par 
\quad \quad \quad $\displaystyle s\leq t\in \rbrack 0,1\lbrack \rightarrow \gamma
 (s):=\int_{{\mathbb{T}}}{\varphi (se^{i\theta })\log ^{-}\left\vert{f(se^{i\theta
 })}\right\vert d\theta }$\par 
is a continuous function of $\displaystyle s\in \lbrack 0,t\rbrack .$
\end{Lmm}
\quad Proof.\ \par 
Because $\displaystyle s\leq t<1,$ the holomorphic function in
 the unit disc $\displaystyle f(se^{i\theta })$ has only a finite
 number of zeroes say $\displaystyle N(t).$ As usual we can factor
 out the zeros of $f$ to get\ \par 
\quad \quad \quad $\displaystyle f(z)=\prod_{j=1}^{N}{(z-a_{j})}g(z)$\ \par 
where $\displaystyle g(z)$ has no zeros in the disc $\displaystyle
 \bar D(0,t).$ Hence we get\ \par 
\quad \quad \quad $\displaystyle \log \left\vert{f(z)}\right\vert =\sum_{j=1}^{N}{\log
 \left\vert{z-a_{j}}\right\vert }+\log \left\vert{g(z)}\right\vert .$\ \par 
Let $\displaystyle a_{j}=r_{j}e^{\alpha _{j}},\ r_{j}>0$ because
 $\displaystyle \ \left\vert{f(0)}\right\vert =1,$ then it suffices
 to show that\ \par 
\quad \quad \quad $\displaystyle \gamma (s):=\int_{{\mathbb{T}}}{\varphi (se^{i\theta
 })\log ^{-}\left\vert{se^{i\theta }-re^{i\alpha }}\right\vert d\theta }$\ \par 
is continuous in $s$ near $\displaystyle s=r,$ because $\displaystyle
 \ \int_{{\mathbb{T}}}{\varphi (se^{i\theta })\log ^{-}\left\vert{g(se^{i\theta
 })}\right\vert d\theta }$ is clearly continuous.\ \par 
\quad To see that $\gamma (s)$ is continuous at $\displaystyle s=r,$
 it suffices to show\ \par 
\quad \quad \quad $\displaystyle \gamma (s_{n})\rightarrow \gamma (r)$ when $\displaystyle
 s_{n}\rightarrow r.$\ \par 
But\ \par 
\quad \quad \quad $\displaystyle \forall \theta \neq 0,\ \varphi (se^{i\theta })\log
 \left\vert{se^{i\theta }-r}\right\vert \rightarrow \varphi (re^{i\theta
 })\log \left\vert{re^{i\theta }-r}\right\vert $\ \par 
and $\displaystyle \log \frac{1}{\left\vert{se^{i\theta }-r}\right\vert
 }\leq c_{\epsilon }\left\vert{se^{i\theta }-r}\right\vert ^{-\epsilon
 }$ with $\displaystyle \epsilon >0.$ So choosing $\displaystyle
 \epsilon <1,$ we get that $\displaystyle \log \frac{1}{\left\vert{se^{i\theta
 }-r}\right\vert }\in L^{1}({\mathbb{T}})$ uniformly in $s.$
 Because $\displaystyle \varphi (se^{i\theta })$ is continuous
 uniformly in $\displaystyle s\in \lbrack 0,t\rbrack $ we get
 also $\displaystyle \varphi (se^{i\theta })\log \frac{1}{\left\vert{se^{i\theta
 }-r}\right\vert }\in L^{1}({\mathbb{T}})$ uniformly in $s.$
 So we can apply the dominated convergence theorem of Lebesgue
 to get the result. $\hfill\blacksquare $\ \par 
\ \par 
\ \par 

\begin{Lmm}
~\label{0_NI2}Suppose that $\displaystyle g_{s}(z)\in {\mathcal{C}}^{\infty
 }(\bar {\mathbb{D}})$ and $\displaystyle f\in {\mathcal{H}}({\mathbb{D}})$
 then, with $\displaystyle s<1,\ f_{s}(z):=f(sz),$ we have:\par 
\quad \quad \quad $\displaystyle \ \sum_{a\in Z(f_{s})}{g_{s}(a)}=\int_{{\mathbb{D}}}{\log
 \left\vert{f_{s}(z)}\right\vert \triangle g_{s}(z)}+\int_{{\mathbb{T}}}{(g_{s}\partial
 _{n}\log \left\vert{f_{s}(z)}\right\vert -\log \left\vert{f_{s}(z)}\right\vert
 \partial _{n}g_{s})}.$
\end{Lmm}
\quad Proof.\ \par 
To apply the Green formula we need $\displaystyle {\mathcal{C}}^{2}(\bar
 {\mathbb{D}})$ functions, so we shall use an approximation of
 $\displaystyle \log \left\vert{f_{s}(z)}\right\vert .$ First
 because $\displaystyle s<1,$ we have that $\displaystyle f_{s}$
 has a finite number of zeroes in $\displaystyle {\mathbb{D}}$
 and we take an $\epsilon >0$ small enough to have the discs
 $\displaystyle \forall a\in Z(f_{s}),\ D(a,\epsilon )$ disjoint.
 Then we consider\ \par 
\quad \quad \quad $\displaystyle u_{\epsilon }(z):=\log \left\vert{f_{s}(z)}\right\vert
 (1-\sum_{a\in Z(f_{s})}{\chi _{a}(z,\epsilon )}),$\ \par 
with $\displaystyle \chi _{a}(z,\epsilon ):=0$ for $\displaystyle
 z\notin D(a,\epsilon ),\ \chi _{a}(z,\epsilon )=1$ for $\displaystyle
 z\in D(a,\epsilon /2),\ 0\leq \chi _{a}(z,\epsilon )\leq 1$
 and $\displaystyle \chi _{a}(z,\epsilon )\in {\mathcal{C}}^{\infty
 }(\bar {\mathbb{D}}).$\ \par 
\quad Then, because $\displaystyle Z(f_{s})$ is finite, we have that
 $\displaystyle u_{\epsilon }$ is in $\displaystyle {\mathcal{C}}^{\infty
 }(\bar {\mathbb{D}})$ and we can apply the Green formula to
 $\displaystyle g_{s}$ and $\displaystyle u_{\epsilon }.$ we have\ \par 
\quad \quad \quad $\displaystyle \ \int_{{\mathbb{D}}}{(g_{s}(z)\triangle u_{\epsilon
 }(z)-u_{\epsilon }(z)\triangle g_{s}(z))}=\int_{{\mathbb{T}}}{(g_{s}(e^{i\theta
 })\partial _{n}u_{\epsilon }(e^{i\theta })-u_{\epsilon }(e^{i\theta
 })\partial _{n}g_{s}(e^{i\theta }))}.$\ \par 
Clearly $\Delta u_{\epsilon }=0$ outside $\displaystyle \ \bigcup_{a\in
 Z(f_{s})}{D(a,\epsilon )}$ and in $\displaystyle D(a,\epsilon
 )$ we get, because $\displaystyle g_{s}(z)$ is continuous in
 $\displaystyle \bar {\mathbb{D}},$\ \par 
\quad \quad \quad $\displaystyle \ \int_{D(a,\epsilon )}{g_{s}(z)\triangle u_{\epsilon
 }(z)}\underset{\epsilon \rightarrow 0}{\rightarrow }g_{s}(a).$\ \par 
We have also\ \par 
\quad \quad \quad $\displaystyle \ \int_{{\mathbb{D}}}{u_{\epsilon }(z)\triangle
 g_{s}(z)}\underset{\epsilon \rightarrow 0}{\rightarrow }\int_{{\mathbb{D}}}{\log
 \left\vert{f_{s}(z)}\right\vert \triangle g_{s}(z)},$\ \par 
\quad \quad \quad $\displaystyle \ \int_{{\mathbb{T}}}{u_{\epsilon }(e^{i\theta
 })\partial _{n}g_{s}(e^{i\theta })}\underset{\epsilon \rightarrow
 0}{\rightarrow }\int_{{\mathbb{T}}}{\log \left\vert{f_{s}(z)}\right\vert
 \partial _{n}g_{s}},$\ \par 
and\ \par 
\quad \quad \quad $\displaystyle \ \int_{{\mathbb{T}}}{g_{s}(e^{i\theta })\partial
 _{n}u_{\epsilon }(e^{i\theta })}\underset{\epsilon \rightarrow
 0}{\rightarrow }\int_{{\mathbb{T}}}{(g_{s}\partial _{n}\log
 \left\vert{f_{s}(z)}\right\vert },$\ \par 
which prove the lemma. $\hfill\blacksquare $\ \par 

\begin{Lmm}
~\label{0_NI3}Suppose that $\displaystyle g_{s}(z)\in {\mathcal{C}}^{\infty
 }(\bar {\mathbb{D}})$ and $\displaystyle u$ is a subharmonic
 function in the disc $\displaystyle {\mathbb{D}}$ ; then, with
 $\displaystyle \forall s<1,\ u_{s}(z):=u(sz),$ we have:\par 
\quad \quad \quad $\displaystyle \ \int_{{\mathbb{D}}}{g_{s}(z)d\mu (z)}=\int_{{\mathbb{D}}}{u_{s}(z)\triangle
 g_{s}(z)}+\int_{{\mathbb{T}}}{(g_{s}\partial _{n}u_{s}-u_{s}\partial
 _{n}g_{s})},$\par 
where $\mu _{s}:=\Delta u_{s}$ is the positive Riesz measure
 associated to $\displaystyle u_{s}.$
\end{Lmm}
\quad Proof.\ \par 
First recall that $\displaystyle \mu :=\Delta u,$ the Riesz measure
 associated to the subharmonic non trivial function $\displaystyle
 u$ in the disc $\displaystyle {\mathbb{D}},$ is finite on the
 compact sets of $\displaystyle {\mathbb{D}}$ because $\displaystyle
 u\in L^{1}_{loc}({\mathbb{D}})$ implies that $\displaystyle
 u\in {\mathcal{D}}'({\mathbb{D}})$ hence $\displaystyle \Delta
 u\in {\mathcal{D}}'({\mathbb{D}})$ ; so take a function $\displaystyle
 \varphi \in {\mathcal{D}}({\mathbb{D}})$ which is $1$ on the
 compact $\displaystyle K\Subset {\mathbb{D}}$ and $\displaystyle
 \varphi \geq 0.$ Then, because $\Delta u$ is a positive measure, we get\ \par 
\quad \quad \quad $\displaystyle \ {\left\langle{\Delta u,\varphi }\right\rangle}=\int_{{\mathbb{D}}}{\varphi
 (z)d\mu (z)}\geq \int_{K}{d\mu (z)}$\ \par 
hence $\displaystyle \mu (K)\leq {\left\langle{\Delta u,\varphi
 }\right\rangle}<\infty .$\ \par 
\quad The idea is to start with the measure $\mu :=\Delta u$ and, because
 $\displaystyle s<1,$ we can cut it by a smooth function $\displaystyle
 \gamma _{s}(z)\in {\mathcal{C}}^{\infty }_{c}({\mathbb{D}}),$
 such that $\gamma (z)=1$ in $\displaystyle D(0,s).$ Then we
 regularise $\gamma \mu $ by convolution with :\ \par 
\quad \quad \quad $\displaystyle \chi _{\epsilon }(\rho e^{i\theta }):=a_{\epsilon
 }(\rho )b_{\epsilon }(\theta ),$\ \par 
with\ \par 
\quad \quad \quad $\displaystyle a_{\epsilon }(\rho ):=\frac{1}{\epsilon }a(\frac{\left\vert{\rho
 }\right\vert }{\epsilon }),\ 0\leq a(t)\leq 1,\ a\in {\mathcal{C}}^{\infty
 }_{c}(\lbrack 0,1\lbrack ),\ t\leq 1/2\Rightarrow a(t)=1.$\ \par 
And\ \par 
\quad \quad \quad $\displaystyle b_{\epsilon }(\rho ):=\frac{1}{\epsilon }a(\frac{\left\vert{\theta
 }\right\vert }{\epsilon }),\ 0\leq b(t)\leq 1,\ a\in {\mathcal{C}}^{\infty
 }_{c}(\lbrack 0,2\pi \lbrack ),\ t\leq 1/2\Rightarrow b(t)=1.$\ \par 
So we set the potential:\ \par 
\quad \quad \quad $\displaystyle U(z):=\int_{{\mathbb{D}}}{\log \left\vert{z-\zeta
 }\right\vert \gamma d\mu (\zeta )}=\log \left\vert{\cdot }\right\vert
 \ast (\gamma \mu )$\ \par 
and we have $\Delta U(z)=\gamma (z)\mu (z)$ in distributions
 sense, and we regularise\ \par 
\quad \quad \quad $\displaystyle U_{\epsilon }:=\chi _{\epsilon }\ast U\Rightarrow
 \Delta U_{\epsilon }=\chi _{\epsilon }\ast \Delta U.$\ \par 
Now we have that $\Delta (u-U)=\mu -\gamma \mu =0$ in $\displaystyle
 D(0,s)$ so $\displaystyle H:=u-U$ is harmonic in $\displaystyle
 D(0,s)$ hence smooth.\ \par 
On the other hand we have, because $\displaystyle U_{\epsilon
 }$ is $\displaystyle {\mathcal{C}}^{\infty },$ that the Green
 formula is applicable so\ \par 
\quad \quad \quad $\displaystyle \ \int_{{\mathbb{D}}}{(g_{s}(z)\triangle U_{\epsilon
 }(sz)-U_{\epsilon }(sz)\triangle g_{s}(z))}=\int_{{\mathbb{T}}}{(g_{s}(e^{i\theta
 })\partial _{n}U_{\epsilon }(se^{i\theta })-U_{\epsilon }(se^{i\theta
 })\partial _{n}g_{s}(e^{i\theta }))}.$\ \par 
And from $\displaystyle u=U+H,$ we get $\displaystyle u=H+\lim
 _{\epsilon \rightarrow 0}U_{\epsilon }$ so it remains to see
 what happen to each term.\ \par 
\quad For the first one\ \par 
\quad \quad \quad $\displaystyle \ \int_{{\mathbb{D}}}{g_{s}(z)\triangle U_{\epsilon
 }(sz)}=\int_{{\mathbb{D}}}{g_{s}(z)(\chi _{\epsilon }\ast \triangle
 U)(sz)}=\int_{{\mathbb{D}}}{(g_{s}\ast \chi _{\epsilon })(\zeta
 )\triangle U)(s\zeta )}.$\ \par 
But $\displaystyle (g_{s}\ast \chi _{\epsilon })(\zeta )\rightarrow
 g_{s}(\zeta )$ uniformly in $\displaystyle \bar {\mathbb{D}},$
 because $\displaystyle g_{s}$ is smooth on $\displaystyle \bar
 {\mathbb{D}},$ and $\displaystyle \Delta U(sz)=\gamma (sz)\mu
 (sz)=\mu (sz)$ is a bounded measure in $\displaystyle \bar {\mathbb{D}}$
 so we get\ \par 
\quad \quad \quad $\displaystyle \ \int_{{\mathbb{D}}}{g_{s}(z)\triangle U_{\epsilon
 }(sz)}\underset{\epsilon \rightarrow 0}{\rightarrow }\int_{{\mathbb{D}}}{g_{s}(z)\triangle
 U_{\epsilon }(sz)}=\int_{{\mathbb{D}}}{g_{s}(z)d\mu (sz)}.$\ \par 
\quad For the second one:\ \par 
\quad \quad \quad $\displaystyle \ \int_{{\mathbb{D}}}{U_{\epsilon }(sz)\triangle
 g_{s}(z)}=\int_{{\mathbb{D}}}{U(s\zeta )(\triangle g_{s}\ast
 \chi _{\epsilon })}\underset{\epsilon \rightarrow 0}{\rightarrow
 }\int_{{\mathbb{D}}}{U(sz)\triangle g_{s}(z)},$\ \par 
as above because $\displaystyle (\Delta g_{s}\ast \chi _{\epsilon
 })(\zeta )\rightarrow \Delta g_{s}(\zeta )$ uniformly in $\displaystyle
 \bar {\mathbb{D}},$ because $\displaystyle \Delta g_{s}$ is
 smooth on $\displaystyle \bar {\mathbb{D}}.$\ \par 
\quad For the third term\ \par 
\quad \quad \quad $\displaystyle \ \int_{{\mathbb{T}}}{g_{s}(e^{i\theta })\partial
 _{n}U_{\epsilon }(se^{i\theta })}=\int_{{\mathbb{T}}}{g_{s}(e^{i\theta
 })(\chi _{\epsilon }\ast \partial _{n}U)(se^{i\theta })}$\ \par 
and here we use the special form of $\chi _{\epsilon }(\rho e^{i\theta
 }):=a_{\epsilon }(\rho )b_{\epsilon }(\theta )$ to get\ \par 
\quad \quad \quad $\displaystyle (\chi _{\epsilon }\ast \partial _{n}U)(sz)=\int_{{\mathbb{D}}}{\chi
 _{\epsilon }(\zeta -z)\partial _{n}U)(\zeta )}=\int_{0}^{1}{a_{\epsilon
 }(\rho -s)\lbrace \int_{{\mathbb{T}}}{b_{\epsilon }(\varphi
 -\theta )\partial _{n}U(\rho e^{i\varphi })d\varphi }\rbrace
 \rho d\rho },$\ \par 
so by Fubini we get\ \par 
\quad \quad \quad $\displaystyle \ \int_{{\mathbb{T}}}{g_{s}(e^{i\theta })(\chi
 _{\epsilon }\ast \partial _{n}U)(se^{i\theta })d\theta }=$\ \par 
\quad \quad \quad \quad \quad \quad \quad $\displaystyle =\int_{{\mathbb{T}}}{\lbrace \int_{{\mathbb{T}}}{g_{s}(e^{i\theta
 })b_{\epsilon }(\varphi -\theta )d\theta \lbrack \int_{0}^{1}{a_{\epsilon
 }(\rho -s)\partial _{n}U(\rho e^{i\varphi })}\rbrack \rho d\rho
 }\rbrace d\varphi }.$\ \par 
But\ \par 
\quad \quad \quad $\displaystyle \ \int_{{\mathbb{T}}}{g_{s}(e^{i\theta })b_{\epsilon
 }(\varphi -\theta )d\theta }=(g_{s}\ast b_{\epsilon })(\varphi )$\ \par 
and \ \par 
\quad \quad \quad $\displaystyle \ \int_{0}^{1}{a_{\epsilon }(\rho -s)\partial
 _{n}U(\rho e^{i\varphi })}\rbrack \rho d\rho \underset{\epsilon
 \rightarrow 0}{\rightarrow }\partial _{n}U(se^{i\varphi })$\ \par 
as a measure on $\displaystyle {\mathbb{T}},$ so, because $\displaystyle
 (g_{s}\ast b_{\epsilon })(\varphi )\underset{\epsilon \rightarrow
 0}{\rightarrow }g_{s}(e^{i\varphi })$ uniformly on $\displaystyle
 {\mathbb{T}}$ because $\displaystyle g_{s}(e^{i\theta })\in
 {\mathcal{C}}^{\infty }({\mathbb{T}}),$ we get\ \par 
\quad \quad \quad $\displaystyle \ \int_{{\mathbb{T}}}{g_{s}(e^{i\theta })(\chi
 _{\epsilon }\ast \partial _{n}U)(se^{i\theta })d\theta }\underset{\epsilon
 \rightarrow 0}{\rightarrow }\int_{{\mathbb{T}}}{g_{s}(e^{i\theta
 })\partial _{n}U)(se^{i\theta })d\theta }.$\ \par 
\quad For the last term, we get the same way:\ \par 
\quad \quad \quad $\displaystyle \ \int_{{\mathbb{T}}}{U_{\epsilon }(se^{i\theta
 })\partial _{n}g_{s}(e^{i\theta })}\underset{\epsilon \rightarrow
 0}{\rightarrow }\int_{{\mathbb{T}}}{\partial _{n}g_{s}(e^{i\theta
 })U(se^{i\theta })d\theta }.$\ \par 
So we get\ \par 
\quad \quad \quad $\displaystyle \ \int_{{\mathbb{D}}}{(g_{s}(z)\triangle U(sz)-U(sz)\triangle
 g_{s}(z))}=\int_{{\mathbb{T}}}{(g_{s}(e^{i\theta })\partial
 _{n}U(se^{i\theta })-U(se^{i\theta })\partial _{n}g_{s}(e^{i\theta
 }))}.$\ \par 
\quad Now we replace $U$ by $\displaystyle U=u+H$ with $H$ harmonic
 in $\displaystyle \bar D(0,s)$ to get\ \par 
\quad \quad \quad $\displaystyle \ \int_{{\mathbb{D}}}{(g_{s}(z)\triangle u(sz)-\lbrack
 u+H\rbrack (sz)\triangle g_{s}(z))}=$\ \par 
\quad \quad \quad \quad \quad \quad \quad $\displaystyle =\int_{{\mathbb{T}}}{(g_{s}(e^{i\theta })\partial
 _{n}\lbrack u+H\rbrack (se^{i\theta })-\lbrack u+H\rbrack (se^{i\theta
 })\partial _{n}g_{s}(e^{i\theta }))},$\ \par 
but because $H$ is $\displaystyle {\mathcal{C}}^{\infty }$ we
 get, applying the Green formula to it\ \par 
\quad \quad \quad $\displaystyle \ \int_{{\mathbb{D}}}{-H(sz)\triangle g_{s}(z))}=\int_{{\mathbb{T}}}{(g_{s}(e^{i\theta
 })\partial _{n}H(se^{i\theta })-H(se^{i\theta })\partial _{n}g_{s}(e^{i\theta
 }))},$\ \par 
so it remains\ \par 
\quad \quad \quad $\displaystyle \ \int_{{\mathbb{D}}}{(g_{s}(z)\triangle u(sz)-u(sz)\triangle
 g_{s}(z))}=\int_{{\mathbb{T}}}{(g_{s}(e^{i\theta })\partial
 _{n}u(se^{i\theta })-u(se^{i\theta })\partial _{n}g_{s}(e^{i\theta
 }))},$\ \par 
which proves the lemma. $\hfill\blacksquare $\ \par 

\begin{Lmm}
~\label{6_A0}Let $\varphi (z)$ be a positive function in $\displaystyle
 {\mathbb{D}}$ and $\displaystyle f\in {\mathcal{H}}({\mathbb{D}})\
 ;$ set $\displaystyle f_{s}(z):=f(sz)$ and suppose that:\par 
\quad \quad \quad $\displaystyle \forall s<1,\ \sum_{a\in Z(f_{s})}{(1-\left\vert{a}\right\vert
 ^{2})^{p+1}\varphi (sa)}\leq \int_{{\mathbb{D}}}{(1-\left\vert{z}\right\vert
 ^{2})^{p-1}\varphi (sz)\log ^{+}\left\vert{f(sz)}\right\vert },$\par 
then, for any $\displaystyle 1>\delta >0$ we have\par 
\quad \quad \quad $\displaystyle \ \sum_{a\in Z(f)}{(1-\left\vert{a}\right\vert
 ^{2})^{p+1}\varphi (a)}\leq \sup _{1-\delta <s<1}\int_{{\mathbb{D}}}{(1-\left\vert{z}\right\vert
 ^{2})^{p-1}\varphi (sz)\log ^{+}\left\vert{f(sz)}\right\vert }.$\par 
\quad We have also:\par 
let $\varphi (z),\ \psi (z)$ be positive continuous functions
 in $\displaystyle {\mathbb{D}}$ and $\displaystyle f\in {\mathcal{H}}({\mathbb{D}})$
 such that:\par 
\quad \quad \quad $\displaystyle \forall s<1,\ \sum_{a\in Z(f)\cap D(0,s)}{(1-\left\vert{a}\right\vert
 ^{2})\varphi (sa)}\leq \int_{{\mathbb{D}}}{\varphi (sz)\log
 ^{+}\left\vert{f(sz)}\right\vert }+\int_{{\mathbb{T}}}{\psi
 (se^{i\theta })\log ^{+}\left\vert{f(se^{i\theta })}\right\vert }$\par 
then, for any $\displaystyle 1>\delta >0$ we have\par 
\quad $\displaystyle \ \sum_{a\in Z(f)}{(1-\left\vert{a}\right\vert
 ^{2})\varphi (a)}\leq \sup _{1-\delta <s<1}\int_{{\mathbb{D}}}{\varphi
 (sz)\log ^{+}\left\vert{f(sz)}\right\vert }+\sup _{1-\delta
 <s<1}\int_{{\mathbb{T}}}{\psi (se^{i\theta })\log ^{+}\left\vert{f(sz)}\right\vert
 }.$
\end{Lmm}
\quad Proof.\ \par 
We have $\displaystyle a\in Z(f_{s})\iff f(sa)=0,$ i.e. $\displaystyle
 b:=sa\in Z(f)\cap D(0,s).$ Hence the hypothesis is\ \par 
\quad \quad \quad $\displaystyle \forall s<1,\ \sum_{a\in Z(f)\cap D(0,s)}{(1-\left\vert{\frac{a}{s}}\right\vert
 ^{2})^{p+1}\varphi (a)}\leq \int_{{\mathbb{D}}}{(1-\left\vert{z}\right\vert
 ^{2})^{p-1}\varphi (sz)\log ^{+}\left\vert{f(sz)}\right\vert }.$\ \par 
We fix $\displaystyle 1-\delta <r<1,\ r<s<1,$ then, because $\displaystyle
 Z(f)\cap D(0,r)\subset Z(f)\cap D(0,s)$ and $\varphi \geq 0,$ we have\ \par 
\quad \quad \quad $\displaystyle \ \sum_{a\in Z(f)\cap D(0,r)}{(1-\left\vert{\frac{a}{s}}\right\vert
 ^{2})^{p+1}\varphi (a)}\leq \sum_{a\in Z(f)\cap D(0,s)}{(1-\left\vert{\frac{a}{s}}\right\vert
 ^{2})^{p+1}\varphi (a)}\leq $\ \par 
\quad \quad \quad \quad \quad \quad \quad \quad \quad \quad \quad \quad \quad \quad \quad \quad \quad \quad \quad $\displaystyle \leq \sup _{1-\delta <s<1}\int_{{\mathbb{D}}}{(1-\left\vert{z}\right\vert
 ^{2})^{p-1}\varphi (z)\log ^{+}\left\vert{f(z)}\right\vert }.$\ \par 
\quad In $\displaystyle D(0,r)$ we have a finite fixed number of zeroes
 of $f,$ and, because $(1-\left\vert{\frac{a}{s}}\right\vert
 ^{2})^{p+1}$ is continuous in $\displaystyle s\leq 1$ for $\displaystyle
 a\in {\mathbb{D}},$ we have\ \par 
\quad \quad \quad $\displaystyle \forall a\in Z(f)\cap D(0,r),\ \lim _{s\rightarrow
 1}(1-\left\vert{\frac{a}{s}}\right\vert ^{2})^{p+1}=(1-\left\vert{a}\right\vert
 ^{2})^{p+1}.$\ \par 
Hence\ \par 
\quad \quad \quad $\displaystyle \ \sum_{a\in Z(f)\cap D(0,r)}{(1-\left\vert{a}\right\vert
 ^{2})^{p+1}\varphi (a)}\leq \sup _{1-\delta <s<1}\int_{{\mathbb{D}}}{(1-\left\vert{z}\right\vert
 ^{2})^{p-1}\varphi (sz)\log ^{+}\left\vert{f(sz)}\right\vert }.$\ \par 
Because the right hand side is independent of $r<1$ and $\varphi
 $ is positive in $\displaystyle {\mathbb{D}}$ so the sequence\ \par 
\quad \quad \quad $\displaystyle S(r):=\sum_{a\in Z(f)\cap D(0,r)}{(1-\left\vert{a}\right\vert
 ^{2})^{p+1}\varphi (a)}$\ \par 
is increasing with $r,$ we get\ \par 
\quad \quad \quad $\displaystyle \ \sum_{a\in Z(f)}{(1-\left\vert{a}\right\vert
 ^{2})^{p+1}\varphi (a)}\leq \sup _{1-\delta <s<1}\int_{{\mathbb{D}}}{(1-\left\vert{z}\right\vert
 ^{2})^{p-1}\varphi (sz)\log ^{+}\left\vert{f(sz)}\right\vert }.$\ \par 
This proves the first part. The proof of the second one is just
 identical. $\hfill\blacksquare $\ \par 
\ \par 

\bibliographystyle{/usr/local/texlive/2013/texmf-dist/bibtex/bst/base/plain}

\end{document}